\documentclass{birkmult}
 \newtheorem{theorem}{Theorem}[section]

 \newtheorem{proposition}[theorem]{Proposition}
 \theoremstyle{definition}
 \newtheorem{definition}[theorem]{Definition}
 \newtheorem{property}{Property}
 \theoremstyle{remark}

 \numberwithin{equation}{section}
 
\usepackage[all]{xy}
\usepackage{amssymb}
\newcommand\cusp{\operatorname{cu}}

\newcommand\sus{\operatorname{sus}}
\newcommand\susn[1]{\operatorname{s}(#1)}
\newcommand\psusn[1]{\operatorname{ps}(#1)}
\newcommand\suss[1]{\operatorname{s\star}(#1)}
\newcommand\tsuss[1]{\operatorname{s\tstar}(#1)}
\newcommand\tstar{\tilde\star}
\newcommand\talpha{\tilde\alpha}
\newcommand\cs{\operatorname{cs}}
\newcommand\cps{\operatorname{cps}}
\newcommand\cpsn[1]{\operatorname{cps}(#1)}
\newcommand\css{\operatorname{cs\star}}
\newcommand\csn[1]{\operatorname{cs}(#1)}

\newcommand\AS{\overline{\operatorname{AS}}}

\newcommand\famtot{M}
\newcommand\fambas{B}

\newcommand\famrel{\famtot/\fambas}

\newcommand\indAS{\operatorname{ind_{\operatorname{AS}}}}

\newcommand\Kto{\operatorname{K^1}}

\newcommand\boxb[1]{\square_b}

\renewcommand{\theequation}{\arabic{equation}}
\let\Appendix\appendix
\newcommand\paperbody%
        {\renewcommand{\theequation}{\thesection.\arabic{equation}}}
\renewcommand\appendix{\Appendix%
\renewcommand{\theequation}{A.\arabic{equation}}
\renewcommand{\thelemma}{A.\arabic{lemma}}}

\newcommand\Kc{K_{\text{c}}}
\newcommand\Kco{K_{\text{c}}^0}

\newcommand\cu{\operatorname{cu}}

\newcommand\cFTs{{}^{\Phi}\overline{T}\kern-1pt{}^*}

\newcommand\Ell{\operatorname{Ell}}

\newcommand\Tr{\operatorname{Tr}}

\newcommand\bTr{\overline{\operatorname{Tr}}}
\newcommand\rdsTr{\operatorname{\overline{Tr}_{\susn{2}}}}
\newcommand\etac{\eta_{\cusp}}

\newcommand\RTr[1]{\operatorname{Tr_{R,#1}}}

\newcommand\Cl{\operatorname{Cl}}

\newcommand\Det{\operatorname{Det}}

\newcommand\tI{\tilde{I}}

\newcommand\com[1]{\overline{#1}}

\newcommand\ie{i\@.e\@. }

\newcommand\cA{\mathcal{A}}

\newcommand\bbC{\mathbb C}

\newcommand\bbN{\mathbb N}

\newcommand\bbR{\mathbb R}
\newcommand\bbS{\mathbb S}

\newcommand\bbZ{\mathbb Z}

\newcommand\CI{{\mathcal{C}}^{\infty}}

\newcommand\cFNs{{}^{\Phi}\overline N\kern-1pt{}^*}

\newcommand\ind{\operatorname{ind}}

\newcommand\tr{\operatorname{tr}}

\newcommand\Hom{\operatorname{Hom}}
\newcommand\Aut{\operatorname{Aut}}
\newcommand\Id{\operatorname{Id}}

\newcommand\pa{\partial}

\renewcommand\Re{\operatorname{Re}}

\newcommand\Mand{\text{ and }}

\newcommand\Mfor{\text{ for }}

\newcommand\Mwhen{\text{ when }}

\newcommand\cf{cf\@. }

\begin{document}
%
%
%
%
%
%
%
%
%
\title[Boundaries, eta invariant and the determinant bundle]
 {Boundaries, eta invariant and the determinant bundle}
\author[Melrose]{Richard Melrose}

\address{%
Department of Mathematics\\
MIT\\
77 Massachusetts Avenue \\ 
Cambridge, MA 02139-4307\\
USA}

\email{rbm@math.mit.edu}

\thanks{The first author acknowledges the support of the National Science
  Foundation under grant DMS0408993,the second author acknowledges support
  of the Fonds qu\'{e}b\'{e}cois sur la nature et les technologies and
  NSERC while part of this work was conducted.}
\author[Rochon]{Fr\'ed\'eric Rochon}
\address{%
Department of Mathematics \\
University of Toronto \\
40 St. George Street  \\
Toronto, Ontario, M5S 2E4 \\
Canada}
\email{rochon@math.utoronto.ca}


\date{September 6 2007}


\begin{abstract}
Cobordism invariance shows that the index, in K-theory, of
  a family of pseudodifferential operators on the boundary of a fibration
  vanishes if the symbol family extends to be elliptic across the whole
  fibration. For Dirac operators with spectral boundary condition, Dai and
  Freed \cite{dai-freed1} gave an explicit version of this at the level
  of the determinant bundle. Their result, that  the eta invariant of the
  interior family trivializes the determinant bundle of the boundary
  family, is extended here to the wider context of pseudodifferential
  families of cusp type.
\end{abstract}

\maketitle

\section*{Introduction}

For a fibration of compact manifolds $\famtot\longrightarrow \fambas$ where
the fibre is a compact manifold with boundary, the cobordism of the index
can be interpreted as the vanishing of the (suspended) family index for
the boundary 
\begin{equation}
\indAS:\Kco(T^*(\pa\famrel))\longrightarrow\Kto(\fambas)
\label{fipomb2.22}\end{equation}
on the image of the restriction map 
\begin{equation}
\Kc(T^*(\famrel))\longrightarrow\Kc(\bbR\times T^*(\pa\famrel))=
\Kco(T^*(\pa\famrel)).
\label{fipomb2.23}\end{equation}
This was realized analytically in \cite{fipomb} in terms of cusp
pseudodifferential operators, namely that any elliptic family of cusp
pseudodifferential operators can be perturbed by a family of order
$-\infty$ to be invertible; this is described as the universal case in
\cite{fresmb1}. For the odd version of \eqref{fipomb2.22} and 
in the special case of Dirac operators, Dai and Freed in \cite{dai-freed1},
showed that the $\tau$ (\ie exponentiated $\eta$) invariant of a
self-adjoint Dirac operator on an odd dimensional compact oriented manifold
with boundary, with augmented Atiyah-Patodi-Singer boundary condition,
defines an element of the inverse determinant line for the boundary Dirac
operator. Here we give a full pseudodifferential version of this, showing
that the $\tau$ invariant for a suspended (hence `odd') family of elliptic
cusp pseudodifferential operators, $P,$ trivializes the determinant
bundle for the indicial family $I(P)$
\begin{equation}
\tau=\exp(i\pi\eta(P)):\Det(I(P)\longrightarrow \bbC^*,
\label{fipomb2.24}\end{equation}
which in this case is a doubly-suspended family of elliptic pseudodifferential
operators; the relation to the Dirac case is discussed in detail.

This paper depends substantially on \cite{perdet} where the determinant on
$2n$ times suspended smoothing families is discussed. This determinant in
the doubly suspended case is used to define the determinant bundle for any
doubly-suspended elliptic family of pseudodifferential operators on a
fibration (without boundary). As in the unsuspended case (see Bismut and
Freed \cite{Bismut-Freed1}), the first Chern class of the determinant
bundle is the 2-form part of the Chern character of the index bundle. The
realization of the eta invariant for singly suspended invertible families
in \cite{MR96h:58169} is extended here to the case of invertible families
of suspended cusp operators. In the Dirac case this is shown to reduce to
the eta invariant of Atiyah, Patodi and Singer for a self-adjoint Dirac
operator with augmented APS boundary condition.

In the main body of the paper the consideration of a self-adjoint Dirac
operator, $\eth,$ is replaced by that of the suspended family, generalizing
$\eth+it,$ where $t$ is the suspension variable. This effectively replaces
the self-adjoint Fredholm operators, as a classifying space for odd
K-theory, by the loop group of the small unitary group (see \cite{APS3},
p.81). One advantage of using suspended operators in this way is that the
regularization techniques of \cite{MR96h:58169} can be applied to define
the eta invariant as an extension of the index. In order to discuss
self-adjoint (cusp) pseudodifferential operators using this suspension
approach, it is necessary to consider somewhat less regular (product-type)
families, generalizing $A+it,$ so we show how to extend the analysis to
this larger setting.

As in \cite{perdet}, we introduce the various determinant bundles in a
direct global form, as associated bundles to principal bundles (of
invertible perturbations) instead of using the original spectral definition
of Quillen \cite{Quillen}. In this way, the fact that the $\tau$ invariant
gives a trivialization of the determinant bundle follows rather directly
from the log-multiplicative property
\begin{equation*}
                  \eta(A*B)=\eta(A)+\eta(B)
\label{fipomb2.96}\end{equation*}
of the eta invariant.

The paper is organized as follows. In Section~\ref{cs.0}, we review the
main properties of cusp operators. In Section~\ref{ld.0}, we consider a
conceptually simpler situation which can be thought as an `even'
counterpart of our result. In Section~\ref{detb}, we present the
determinant bundle as an associated bundle to a principal bundle; this
definition is extended to family of $2n$-suspended elliptic operators in
Section~\ref{2n-sus}. This allows us in Section~\ref{totdlb} to rederive a
well-known consequence of the cobordism invariance of the index at the
level of determinant bundles using the contractibility result of
\cite{fipomb}. In Section~\ref{tsc.0}, we introduce the notion of cusp
suspended $*$-algebra, which is used in Section~\ref{Lifting} to lift the
determinant from the boundary. This lifted determinant is defined using the
eta invariant for invertible suspended cusp operators introduced in
Section~\ref{eta-inv}. In Section~\ref{Tdb.0}, we prove the trivialization
result and in Section~\ref{Diracfami} we relate it to the result of Dai and
Freed \cite{dai-freed1} for Dirac operators.  Finally, in
Section~\ref{asg.0}, these results are extended to include the case of a
self-adjoint family of elliptic cusp pseudodifferential operators. This
involves the use of product-type suspended operators, which are discussed
in the Appendix.

\paperbody
\section{Cusp pseudodifferential operators} \label{cs.0}

This section is intended to be a quick summary of the main properties of
cusp pseudodifferential operators and ellipticity. We refer to
\cite{Mazzeo-Melrose4}, \cite{Melrose-Nistor} and \cite{fipomb} for more
details.

Let $Z$ be a compact manifold with non-empty boundary $\pa Z.$  Let 
$x\in\CI(Z)$ be a defining function for the boundary, that is, $x\ge 0$
everywhere on $Z,$ 
\[
\pa Z=\{z\in Z;x(z)=0\}
\]
and $dx(z)\ne 0$ for all $z\in
\pa $Z. Such a choice of boundary defining function determines a cusp
structure on the manifold $Z,$ which is an identification of the normal bundle
of the boundary $\pa Z$ in $Z$ with $\pa Z\times L$ for a 1-dimensional
real vector space $L.$ If $E$ and $F$ are complex vector bundles on $Z,$
then $\Psi_{\cusp}^{m}(Z;E,F)$ denotes the space of cusp pseudodifferential
operators acting from $\CI(Z;E)$ to $\CI(Z;F)$ associated to the choice of
cusp structure. Different choices lead to different algebras
of cusp pseudodifferential operators, but all are isomorphic. We therefore
generally ignore the particular choice of cusp structure.

A \emph{cusp vector field} $V\in \CI(Z,TZ)$ is a vector field such that
$Vx\in x^{2}\CI(Z)$ for any defining function consistent with the chosen
cusp structure. We denote by $\mathcal{V}_{\cusp}(Z)$ the Lie algebra of
such vector fields. The \emph{cusp tangent bundle} ${}^{\cusp}TZ$ is the
smooth vector bundle on $Z$ such that
$\mathcal{V}_{\cusp}=\CI(Z;{}^{\cusp}TZ);$ it is isomorphic to $TZ$ as a
vector bundle, but not naturally so.

Let ${}^{\cusp}S^{*}Z=({}^{\cusp}T^{*}Z\setminus 0)/\bbR^{+}$ be the
\emph{cusp cosphere bundle} and let $R^{m}$ be the trivial complex line bundle 
on ${}^{\cusp}S^{*}Z$ with sections given by functions over ${}^{\cusp}T^{*}Z
\setminus 0$ which are positively homogeneous of degree $m.$

\begin{proposition}\label{cusp.3}
For each $m\in \bbZ$, there is a \emph{symbol map} giving a short exact sequence
\begin{equation}
\Psi_{\cusp}^{m-1}(Z;E,F)
\longrightarrow \Psi_{\cusp}^{m}(Z;E,F)\overset{\sigma_m}
\longrightarrow \CI({}^{\cusp}S^{*}Z;\hom(E,F)\otimes R^{m}).
\label{cusp.4}\end{equation}
\end{proposition}

Then $A\in\Psi_{\cusp}^{m}(Z;E,F)$ is said to be \emph{elliptic} if its
symbol is invertible. In this context, ellipticity is not a sufficient
condition for an operator of order $0$ to be Fredholm on $L^2.$

More generally, one can consider the space of full symbols of order $m$ 
\begin{equation*}
\mathcal{S}_{\cusp}^{m}(Z;E,F)= \rho^{-m}\CI(\overline{{}^{\cusp}T^{*}Z};
\hom(E,F)), 
\label{fipomb2.97}\end{equation*}
where $\rho$ is a defining function for the boundary (at infinity) in the radial 
compactification of ${}^{\cusp}T^{*}Z.$ After choosing appropriate metrics
and connections, one can define a quantization map following standard
constructions
\begin{equation}
 q: \mathcal{S}_{\cusp}^{m}(Z;E,F)\longrightarrow \Psi^{m}_{\cusp}(Z;E,F)
\label{qm.1}\end{equation}
which induces an isomorphism of vector spaces
\begin{equation}
       \mathcal{S}_{\cusp}^{m}(Z;E,F)/\mathcal{S}_{\cusp}^{-\infty}(Z;E,F)\cong
           \Psi_{\cusp}^{m}(Z;E,F)/\Psi_{\cusp}^{-\infty}(Z;E,F).
\label{qm.1p}\end{equation}

If $Y$ is a compact manifold without boundary and $E$ is a complex vector
bundle over $Y,$ there is a naturally defined algebra of suspended
pseudodifferential operators, which is denoted here
$\Psi^{*}_{\sus}(Y;E).$ For a detailed discussion of this algebra (and the
associated modules of operators between bundles) see \cite{MR96h:58169}. An
element $A\in\Psi_{\sus}^{m}(Y;E)$ is a one-parameter family of
pseudodifferential operators in $\Psi^{m}(Y;E)$ in which the parameter
enters symbolically. A suspended pseudodifferential operator is associated
to each cusp pseudodifferential operator by `freezing coefficients at the
boundary.' Given $A\in\Psi_{\cusp}^{m}(Z;E,F)$, for each
$u\in\CI(Z;E)$, $Au\big|_{\pa Z}\in\CI(\pa Z;F)$
depends only on $u\big|_{\pa Z}\in\CI(\pa Z;E).$ The resulting operator
$A_{\pa}:\CI(\pa Z;E)\longrightarrow \CI(\pa Z;F)$ is an element of
$\Psi^m(\pa Z;E,F).$ More generally, if $\tau\in\bbR$ then
\begin{multline}
\Psi_{\cusp}^{m}(Z;E,F)\ni A\longmapsto
e^{i\frac{\tau}{x}}Ae^{-i\frac{\tau}{x}}\in \Psi_{\cusp}^{m}(Z;E,F)\Mand\\
I(A,\tau)=(e^{i\frac{\tau}{x}}Ae^{-i\frac{\tau}{x}})_{\pa}\in
\Psi_{\sus}^{m}(\pa Z;E,F)
\label{cusp.7}\end{multline}
is the \emph{indicial family} of $A.$

\begin{proposition}\label{cusp.11}The \emph{indicial homomorphism} gives a
short exact sequence, 
\begin{equation*}
x\Psi_{\cusp}^{m}(Z;E,F)\longrightarrow \Psi_{\cusp}^{m}(Z;E,F)\overset{I}
\longrightarrow \Psi_{\sus}^{m}(\pa Z; E,F).
\label{fipomb2.52}\end{equation*}
\end{proposition} 

There is a power series expansion for operators $A\in\Psi_{\cusp}^{m}(Z;E,F)$
at the boundary of which $I(A)$ is the first term. Namely, if $x$ is a
boundary defining function consistent with the chosen cusp structure there
is a choice of product decomposition near the boundary consistent with $x$
and a choice of identifications of $E$ and $F$ with their restrictions to
the boundary. Given such a choice the `asymptotically
translation-invariant' elements of $\Psi_{\cusp}^{m}(Z;E,F)$ are
well-defined by 
\begin{equation}
[x^2D_x,A]\in x^{\infty}\Psi_{\cusp}^{m}(Z;E,F)
\label{fipomb2.42}\end{equation}
where $D_x$ acts through the product decomposition. In fact 
\begin{equation}
\left\{A\in\Psi_{\cusp}^{m}(Z;E,F);\text{ \eqref{fipomb2.42} holds}\right\}/
x^\infty\Psi_{\cusp}^{m}(Z;E,F)\overset{I}\longrightarrow 
\Psi_{\sus}^{m}(\pa Z; E,F)
\label{fipomb2.53}\end{equation}
is an isomorphism. Applying Proposition~\ref{cusp.11} repeatedly and using this
observation, \emph{any} element of $\Psi_{\cusp}^{m}( Z; E,F)$ then has a
power series expansion
\begin{equation}
A\sim\sum\limits_{j=0}^{\infty}x^{j}A_j,\
A_j\in\Psi_{\cusp}^{m}(Z;E,F),\ 
[x^2D_x,A_j]\in x^{\infty}\Psi_{\cusp}^{m}(Z;E,F)
\label{fipomb2.43}\end{equation}
which determines it modulo $x^{\infty}\Psi_{\cusp}^{m}(Z;E,F).$ Setting
$I_j(A)=I(A_j)$ this gives a short exact sequence 
\begin{equation}
\begin{gathered}
\xymatrix@1{
x^\infty\Psi_{\cusp}^{m}(Z;E,F)\ar[r]&
\Psi_{\cusp}^{m}(Z;E,F)\ar[r]^(.43){I_*}&
\Psi_{\sus}^{m}(\pa Z; E,F)[[x]],}\\
I_*(A)=\sum\limits_{j=0}^{\infty}x^jI_j(A)
\end{gathered}
\label{fipomb2.44}\end{equation}
which is multiplicative provided the image modules are given the induced
product
\begin{equation}
I_{*}(A)*I_{*}B=
\sum_{j=0}^{\infty} \frac{(ix^{2})^{j}}{j!} (D_{\tau}^{j}I_{*}(A))
                                          (D_{x}^{j}I_{*}(B)).
\label{star-product}\end{equation}
This is equivalent to a star product although not immediately in the
appropriate form because of the asymmetry inherent in \eqref{fipomb2.44};
forcing the latter to be symmetric by iteratively commuting $x^{j/2}$ to
the right induces an explicit star product in $x^2.$ In contrast to 
Proposition~\ref{cusp.11}, the sequence \eqref{fipomb2.44} does depend on
the choice of product structure, on manifold and bundles, and the choice of
the defining function.

A cusp pseudodifferential operator $A\in \Psi_{\cusp}^{m}(Z;E,F)$ is
said to be \emph{fully elliptic} if it is elliptic and if its indicial
family $I(A)$ is invertible in $\Psi_{\sus}^{*}(\pa Z; E,F);$ this is
equivalent to the invertibility of $I(A,\tau)$ for each $\tau$ and to
$I_*(A)$ with respect to the star product.

\begin{proposition}
A cusp pseudodifferential operator is Fredholm acting on the natural cusp
Sobolev spaces if and only if it is fully elliptic.
\label{cusp.12}\end{proposition}

For bundles on a compact manifold without boundary, let
$G^m_{\sus}(Y;E,F)\subset \Psi^m_{\sus}(Y;E,F)$ denote the subset of
elliptic and invertible elements. The $\eta$ invariant of Atiyah, Patodi
and Singer, after reinterpretation, is extended in \cite{MR96h:58169} to a map
\begin{equation}
\begin{gathered}
\eta:G^m_{\sus}(Y;E,F)\longrightarrow \bbC,\\
\eta (AB)=\eta (A)+\eta (B),\
A\in G^m_{\sus}(Y;F,G),\ B\in G^{m'}_{\sus}(Y;E,F).
\end{gathered}
\label{fipomb2.54}\end{equation}
In \cite{Melrose-Nistor} an index theorem for fully elliptic fibred cusp
operators is obtained, as a generalization of the Atiyah-Patodi-Singer
index theorem.

\begin{theorem} [\cite{Melrose-Nistor}] Let $P \in \Psi_{\cusp}^{m}(X;E,F)$
be a fully elliptic operator, then the index of $P$ is given by the formula
\begin{equation} 
\ind(P)=\AS(P)-\frac12\eta(I(P))
\label{ni.3}\end{equation}
where $\AS$ is a regularized integral involving only a finite number of
terms in the full symbol expansion of $P,$ $I(P)\in\Psi_{\sus}^{m}(\pa
X;E)$ is the indicial family of $P$ and $\eta$ is the functional \eqref{fipomb2.54}
introduced in \cite{MR96h:58169}.
\label{ni.1}\end{theorem}

Note that the ellipticity condition on the symbol of $P$ implies that $E$
and $F$ are isomorphic as bundles over the boundary, since $\sigma _m(P)$
restricted to the inward-pointing normal gives such an isomorphism. Thus
one can freely assume that $E$ and $F$ are identified near the boundary.

In the case of a Dirac operator arising from a product structure near the
boundary with invertible boundary Dirac operator and spectral boundary
condition, the theorem applies by adding a cylindrical end on which the
Dirac operator extends to be translation-invariant, with the indicial
family becoming the spectral family for the boundary Dirac operator (for
pure imaginary values of the spectral parameter). The formula \eqref{ni.3}
then reduces to the Atiyah-Patodi-Singer index theorem.

The result \eqref{ni.3} is not really in final form, since the integral
$\AS(P)$ is not given explicitly nor interpreted in any topological
sense. However, since it is symbolic, $\AS(P)$ makes sense if $P$ is only
elliptic, without assuming the invertibility of the indicial family. It
therefore defines a smooth function
\begin{equation}
\AS:\Ell_{\cusp}^m(X;E,F)\longrightarrow \bbC
\label{fipomb2.3}\end{equation}
for each $m.$ We show in Theorem~\ref{fipomb2.4} below that this function
is a log-determinant for the indicial family.

Cusp operators of order $-\infty$ are in general not compact, so in
particular not of trace class. Nevertheless, it is possible to define a
regularized trace which will be substantially used in this paper. 
\begin{proposition}
For $A\in \Psi_{\cusp}^{-n-1}(Z),$ $n=\dim(Z)$ and $z\in \bbC$, the function
$z\mapsto \Tr(x^{z}A)$ is holomorphic for $\Re z> 1$ and has a meromorphic 
extension to the whole complex plane with at most simple poles at $1-\bbN_{0},$
$\bbN_{0}=\{0,1,2,\ldots\}.$
\label{cusp.13}\end{proposition}

For $A\in \Psi_{\cusp}^{-n-1}(Z),$ the \emph{boundary residue trace} of $A,$
denoted $\RTr\pa (A)$, is the residue at $z=0$ of the meromorphic function
$z\mapsto \Tr(x^{z}A).$ In terms of the expansion \eqref{fipomb2.43}
\begin{equation}
\RTr\pa(A)= \frac{1}{2\pi} \int_{\bbR} \Tr(I_1(A,\tau))d\tau .
\label{fipomb2.57}\end{equation}

The \emph{regularized trace} is defined to be 
\begin{equation*}
     \bTr(A)=\lim_{z\to 0} \left( \Tr(x^{z}A)-\frac{\RTr\pa(A)}{z}\right),\Mfor
A\in \Psi_{\cusp}^{-n-1}(Z).
\label{fipomb2.58}\end{equation*}
For $A\in x^{2}\Psi_{\cusp}^{-n-1}(Z)$ this reduces to the usual trace but in
general it is not a trace, since it does not vanish on all
commutators. Rather, there is a \emph{trace-defect formula}
\begin{multline}
\bTr([A,B])=\frac{1}{2\pi i} \int_{\bbR} \Tr\left(I(A,\tau)\frac{\pa}{\pa \tau}
             I(B,\tau)\right)d\tau,\\ A\in\Psi_{\cusp}^{m}(Z),\
	     B\in\Psi_{\cusp}^{m'}(Z),\ m+m'\le -n-1.
\label{cusp.18}\end{multline}
The sign of this formula is correct provided we use \eqref{cusp.7} to
define the indicial family.  Notice that there is a (harmless) sign mistake
in the trace-defect formula of \cite{fipomb}, where a different convention
for the indicial family is used.

\section{Logarithm of the determinant}\label{ld.0}

As a prelude to the discussion of the determinant bundle, we will consider
the conceptually simpler situation of the principal $\bbZ$-bundle
corresponding to the 1-dimensional part of the odd index.
We first recall the generalization of the notion of principal bundle
introduced in \cite{perdet}.

\begin{definition}
Let $G$ be a smooth group (possibly infinite dimensional), then a 
smooth fibration $\mathcal{G}\longrightarrow B$ over a compact manifold $B$ with 
typical fibre $G$ is called a \textbf{bundle of groups} with model $G$ if
its structure group is contained in $\Aut(G)$, the group of smooth
automorphisms of $G.$
\label{gpb.1}\end{definition}

\begin{definition}
Let $\phi:\mathcal{G}\longrightarrow B$ be a bundle of groups with model $G,$ then a
(right) \textbf{principal $\mathcal{G}$-bundle}
is a smooth fibration $\pi: \mathcal{P}\longrightarrow B$ with typical fibre $G$ 
together with a smooth fibrewise (right) group action 
\begin{equation*}
h:\mathcal{P}\times_{B} \mathcal{G}\ni(p,g) \longmapsto p\cdot g\in\mathcal{P}
\label{fipomb2.99}\end{equation*}
which is free and transitive, where  
\begin{equation*}
\mathcal{P}\times_{B} \mathcal{G}= \{ (p,g)\in \mathcal{P}\times \mathcal{G};
\quad \pi(p)= \phi(g) \}.
\label{fipomb2.100}\end{equation*}
\label{gpb.2}\end{definition}
In particular, a principal $G$-bundle $\pi: \mathcal{P}\longrightarrow B$ is automatically
a principal $\mathcal{G}$-bundle where $\mathcal{G}$ is the 
trivial bundle of groups 
\begin{equation*}
                \mathcal{G}= G\times B\longrightarrow B
\label{fipomb2.101}\end{equation*}
given by the projection on the right factor.  In that sense, 
definition~\ref{gpb.2}  is a generalization of the notion of a principal 
bundle.

Notice also that given a bundle of groups $\mathcal{G}\longrightarrow B$, then
$\mathcal{G}$ itself is a principal $\mathcal{G}$-bundle.  It is the
\textbf{trivial principal $\mathcal{G}$-bundle}.  More generally, we say that 
a principal $\mathcal{G}$-bundle $\mathcal{P}\longrightarrow B$ is 
\textbf{trivial} 
if there exists a diffeomorphism $\Psi: \mathcal{P}\longrightarrow \mathcal{G}$ which
preserves the fibrewise group action: 
\begin{equation*}
            \Psi(h(p,g))= \Psi(p)g, \quad \forall (p,g)\in 
  \mathcal{P}\times_{B}\mathcal{G}.
\label{fipomb2.102}\end{equation*}

In this section, the type of principal $\mathcal{G}$-bundle of interest
arises by considering an elliptic family $Q\in\Psi_{\sus}^{m}(M/B;E,F)$ of
suspended operators over a fibration
\begin{equation}
\xymatrix{Y\ar@{-}[r]&M\ar[d]^{\phi}\\&B}
\label{lodet.1}\end{equation}
of compact manifolds without boundary (not necessarily bounding a fibration 
with boundary). Namely, it is given by the smooth fibration
$\mathcal{Q}\longrightarrow B,$ with fibre at $b\in \fambas$
\begin{multline}
\mathcal{Q}_{b}=\\
\left\{Q_{b}+R_{b};R_{b}\in\Psi^{-\infty}_{\sus}(Y_{b},E_{b},F_{b});\
\exists\ (Q_{b}+R_{b})^{-1}\in\Psi^{-m}_{\sus}(Y_{b};F_{b},E_{b})\right\},
\label{lodet.3}\end{multline} 
the set of all invertible perturbations of $Q_b.$ The fibre is non-empty
and is a principal space for the action of the once-suspended smoothing group 
\begin{equation}
G_{\sus}^{-\infty}(Y_b;E_b)=\left\{\Id+A;A\in\Psi_{\sus}^{-\infty}(Y;E_b),\
(\Id+A)^{-1}\in\Id+\Psi_{\sus}^{-\infty}(Y;E_b)\right\}
\label{fipomb2.25}\end{equation}
acting on the right. Thus, $\mathcal{Q}$ is a principal 
$G_{\sus}^{-\infty}(M/B;E)$-bundle with respect to the bundle of groups
$G_{\sus}^{-\infty}(M/B;E)\longrightarrow B$ with fibre at $b\in B$ given by
\eqref{fipomb2.25}.

The structure group at each point is a classifying space for even K-theory
and carries an index homomorphism
\begin{equation}
\ind:G_{\sus}^{-\infty}(Y_{b};E_b)\longrightarrow \bbZ,\
\ind(\Id+A)=\frac1{2\pi i}\int_{\bbR}
\Tr\left(\frac{dA(t)}{dt}(\Id+A(t))^{-1}\right)dt
\label{fipomb2.26}\end{equation}
labelling the components, \ie giving the 0-dimensional cohomology. For a
suspended elliptic family this induces an integral 1-class on $B;$ namely
the first Chern class of the odd index bundle of the family. This can be
seen in terms of the induced principal $\bbZ$-bundle $\mathcal{Q}_{\bbZ}$
associated to $\mathcal{Q}$ 
\begin{equation}
\mathcal{Q}_{\bbZ}= \mathcal{Q}\times\bbZ / \sim,\quad
(Ag,m-\ind(g))\sim (A,m),\quad  \forall g\in G_{\sus}^{-\infty}(Y;E_b).
\label{fipomb2.27}\end{equation}
Since $\bbC^*=\bbC\setminus\{0\}$ is a classifying space for $\bbZ,$ such
bundles are classified up to equivalence by the integral 1-cohomology of
the base.

More explicitly, any principal $\bbZ$-bundle $\phi:P\longrightarrow B$ admits a
`connection' in the sense of a map $h: P\longrightarrow \bbC$ such that $h(mp)=h(p)+m$
for the action of $m\in \bbZ.$ Then the integral 1-class of the principal
$\bbZ$-bundle $P$ is given by the map
\begin{equation}
e^{2\pi ih}:\fambas\longrightarrow \bbC^*
\label{fipomb2.32}\end{equation}
or the cohomology class of $dh$ seen as a 1-form on $B.$ The triviality of
the principal $\bbZ$-bundle is equivalent to the vanishing of the integral
1-class, that is, to the existence of a function
$f:\fambas\longrightarrow\bbC$ such that $h-\phi^{*}f$ is locally constant. 

Moreover, restricted to the `residual' subgroup $G_{\sus}^{-\infty}(Y;E),$
the eta functional of \eqref{fipomb2.54} reduces to twice the index 
\begin{equation}
    \eta\big|_{G_{\sus}^{-\infty}(Y;E)}=2\ind.
\label{fipomb2.29}\end{equation}

In case the fibration is the boundary of a fibration of compact manifolds
with boundary, as in \cite{fipomb}, and the suspended family is the indicial
family of an elliptic family of cusp pseudodifferential operators then we
know that the whole odd index of the indicial family vanishes in odd
K-theory. In particular the first Chern class vanishes and the associated
principal $\bbZ$-bundle is trivial. 

\begin{theorem}\label{fipomb2.4} The eta invariant defines a connection
$\frac12\eta(A)$ on the principal $\bbZ$-bundle in
\eqref{fipomb2.27} (so the first odd Chern class is $\frac12d\eta)$ and in
the case of the indicial operators of a family of elliptic cusp operators, the
Atiyah-Singer term in the index formula \eqref{ni.3} is a log-determinant
for the indicial family, so trivializing the $\bbZ$-bundle.
\end{theorem}

\begin{proof} By \eqref{fipomb2.29}, the function on $\mathcal{Q}\times\bbZ$ 
\begin{equation}
            h(A, m)= \frac12\eta(A)+ m
\label{fipomb2.30}\end{equation}
descends to $\mathcal{Q}_{\bbZ}$ and defines a connection on it. Thus the map
\begin{equation}
\tau=\exp(i\pi\eta):\fambas\longrightarrow \bbC^{*} 
\label{fipomb2.31}\end{equation}
gives the classifying 1-class, the first odd Chern class in
$H^1(\fambas,\bbZ)$ of the index bundle. In general this class is not
trivial, but when $Q=I(Q_{\cusp})$ is the indicial family of a family of
fully elliptic cusp operators $Q_{\cusp}$, the Atiyah-Singer term
$\AS(Q_{\cusp})$ is a well-defined smooth function which does not
depend on the choice of the indicial family modulo $G_{\sus}^{-\infty}(Y;E)$.
>From  formula \eqref{ni.3}
\begin{equation}
h-\AS(Q_{\cusp})=-\AS(Q_{\cusp})+\frac12\eta(A)+m= -\ind(A_{\cusp},b)+m
\label{lodet.5}\end{equation}
is locally constant. This shows that the Atiyah-Singer term explicitly
trivializes the principal $\bbZ$-bundle $\mathcal{Q}_{\bbZ}.$
\end{proof}

\section{The determinant line bundle}
\label{detb}

Consider a fibration of closed manifolds as in \eqref{lodet.1} and let $E$
and $F$ be complex vector bundles on $M.$ Let $P\in\Psi^{m}(M/B;E,F)$ be a
smooth family of elliptic pseudodifferential operators acting on the
fibres.  If the numerical index of the family vanishes, then one can, for
each $b\in B,$ find $Q_{b}\in \Psi^{-\infty}(Y_{b};E_{b},F_{b})$ such that
$P_{b}+Q_{b}$ is invertible. The families index, which is an element of the 
even K-theory of the base $K^{0}(B)$ (see \cite{Atiyah-Singer4} for a definition), 
is the obstruction to the
existence of a smooth family of such perturbations. This obstruction can be
realized as the non-triviality of the bundle with fibre
\begin{equation}
\mathcal{P}_{b}=\left\{P_{b}+Q_{b};\ Q_{b}\in\Psi^{-\infty}(Y_{b},E_{b},F_{b}),\ 
\exists\ (P_{b}+Q_{b})^{-1}\in\Psi^{-k}(Y_{b};F_{b},E_{b})\right\}.
\label{detb.3}\end{equation} 
As in the odd case discussed above, the fibre is non-trivial (here because the
numerical index is assumed to vanish) and is a bundle of principal $G$-spaces
for the groups
\begin{equation}
G^{-\infty}(Y_{b};E)=\left\{\Id + Q;Q\in \Psi^{-\infty}(Y_{b};E),\
\exists\ (\Id +Q)^{-1}\in \Psi^{0}(Y_{b};E)\right\}
\label{det.5}\end{equation}
acting on the right. Thus, $\mathcal{P}\longrightarrow B$ is a principal
$G^{-\infty}(M/B;E)$-bundle for the bundle of groups
$G^{-\infty}(M/B;E)\longrightarrow B$ with fibre at $b\in B$ given by
\eqref{det.5}.

The Fredholm determinant 
\begin{equation*}
\det:\Id+\Psi^{-\infty}(X;W)\longrightarrow \bbC
\label{fipomb2.33}\end{equation*}
is well-defined for any compact manifold $X$ and vector bundle $W.$ It is
multiplicative  
\begin{equation*}
\det(AB)=\det(A)\det(B)
\label{fipomb2.34}\end{equation*}
and is non-vanishing precisely on the group $G^{-\infty}(X;W).$ Explicitly,
it may be defined by 
\begin{equation}
\det(B)=\exp \left( \int_{[0,1]} \gamma^{*}\Tr(A^{-1}dA)\right)
\label{detb.7}\end{equation} 
where $\gamma:[0,1]\longrightarrow G^{-\infty}(X;W)$ is any smooth path with 
$\gamma(0)=\Id$ and $\gamma(1)=B.$  Such a path exists since 
$G^{-\infty}(X;W)$ is connected and the result does not 
depend on the choice of $\gamma$ in view of the integrality of the 1-form
$\frac1{2\pi i}\Tr(A^{-1}dA)$ (which gives the index for the loop group).

\begin{definition}
If $P\in\Psi^{m}(M/B;E,F)$ is a family of elliptic pseudodifferential
operators with vanishing numerical index and $\mathcal{P}\longrightarrow B$ is the
bundle given by \eqref{detb.3}, then the determinant line bundle 
$\Det(P)\longrightarrow B$ of $P$ is the associated line bundle given by
\begin{equation}
\Det(P)=\mathcal{P}\times_{G^{-\infty}(M/B;E)} \bbC
\label{detb.9}\end{equation}
where $G^{-\infty}(Y_{b};F_{b})$ acts on $\bbC$ via the determinant; thus,
$\Det(P)$ is the space $\mathcal{P}\times\bbC$ with the equivalence
relation 
\begin{equation*}
              (A,c)\sim (Ag^{-1}, \det(g)c)
\label{fipomb2.98}\end{equation*}
for $A\in\mathcal{P}_{b},$ $g\in G^{-\infty}(Y_{b};F_{b}),$ $b\in B$ and 
$c\in \bbC.$ 
\label{detb.8}\end{definition}

\noindent As discussed in \cite{perdet}, this definition is equivalent to
the original spectral definition due to Quillen \cite{Quillen}.

If $P\in\Psi^{m}(M/B;E,F)$ is a general elliptic family, with possibly
non-vanishing numerical index, it is possible to give a similar definition
but depending on some additional choices. Assuming for definiteness that
the numerical index is $l\ge0$ one can choose a trivial $l$-dimensional
subbundle $K\subset\CI(\famrel;E)$ as a bundle over $B$, a Hermitian
inner product on $E$ and a volume form on $B.$ Then the fibre in
\eqref{detb.3} may be replaced by
\begin{equation}
\mathcal{P}_{b,K}=\left\{P_{b}+Q_{b};\
Q_{b}\in\Psi^{-\infty}(Y_{b},E_{b},F_{b}),\, \ker (P_b+Q_b)=K_b\right\}.
\label{fipomb2.35}\end{equation}
This fibre is non-empty and for each such choice of $Q_b$ there is a unique
element $L_b\in\Psi^{-m}(Y_{b};F_{b},E_{b})$ which is a left inverse of
$P_b+Q_b$ with range $K_b^\perp$ at each point of $\fambas.$ The action of
the bundle of groups $G^{-\infty}(\famrel;F)$ on the left makes this into a
(left) principal $G^{-\infty}(M/B;F)$-bundle. Then the fibre of the
determinant bundle may be taken to be 
\begin{equation}
\Det(P)_{b,K}=\mathcal{P}_{b,K}\times\bbC/\sim,\ (A,c)\sim (BA,\det(B)c).
\label{fipomb2.36}\end{equation}
In case the numerical index is negative there is a similar construction
intermediate between the two cases.

\section{The $2n$-suspended determinant bundle}\label{2n-sus}

As described in \cite{perdet}, it is possible to extend the notion of
determinant, and hence that of the determinant line bundle, to suspended
pseudodifferential operators with an even number of parameters.

Let $L\in\Psi_{\susn{2n}}^{m}(M/B;E,F)$ be an elliptic family of
$(2n)$-suspended pseudodifferential operators. Ellipticity (in view of the
symbolic dependence on the parameters) implies that such a family is
invertible near infinity in $\bbR^{2n}.$ Thus the families index is
well-defined as an element of the compactly supported K-theory
$\Kc(\bbR^{2n})=\bbZ.$ By Bott periodicity this index may be
identified with the numerical index of a family where the parameters are
quantized, see the discussion in \cite{perdet}. Even assuming the vanishing
of this numerical index, to get an explicitly defined determinant bundle,
as above, we need to introduce a formal parameter $\epsilon$.  

Let $\Psi^{m}_{\susn{2n}}(Y;F)[[\epsilon]]$ denote the space of formal power 
series in $\epsilon$ with coefficients in $\Psi^m_{\susn{2n}}(Y;F).$ For $A\in
\Psi_{\susn{2n}}^{m}(Y;F)[[\epsilon]]$ and $B\in
\Psi_{\susn{2n}}^{m'}(Y;F)[[\epsilon]],$ consider the \emph{$*$-product} 
$A*B\in\Psi^{m+m'}_{\susn{2n}}(Y;F)[[\epsilon]]$ given by
\begin{equation}
\begin{aligned}
A*B(u) &= (\sum_{\mu=0}^{\infty}a_{\mu}\epsilon^{\mu}) * 
          (\sum_{\nu=0}^{\infty}b_{\nu}\epsilon^{\nu}) \\
       &=\sum_{\mu=0}^{\infty}\sum_{\nu=0}^{\infty}\epsilon^{\mu+\nu}
           \left( \sum_{p=0}^{\infty}\frac{i^{p}\epsilon^{p}}{2^{p}p!}
              \omega(D_{v},D_{w})^{p}A(v)B(w)\big|_{v=w=u}\right)
\end{aligned}
\label{detb.10}\end{equation}
where $u,v,w\in\bbR^{2n}$ and $\omega$ is the standard symplectic form on
$\bbR^{2n},$ $\omega(v,w)=v^{T}Jw$ with
\begin{equation}
J=
\begin{pmatrix}
0 & -\Id_{n} \\
\Id_{n} & 0   
\end{pmatrix}.
\label{symp.0}\end{equation}
That \eqref{detb.10} is an associative product follows from its
identification with the usual `Moyal product' arising as the symbolic
product for pseudodifferential operators on $\bbR^n.$

\begin{definition}\label{detb.9a} The module
  $\Psi^{m}_{\susn{2n}}(Y;E,F)[[\epsilon]]$ with $*$-product as in 
  \eqref{detb.10} will be denoted
  $\Psi^{m}_{\suss{2n}}(Y;E,F)[[\epsilon]]$ and the quotient by the
  ideal $\epsilon^{n+1}\Psi^{m}_{\susn{2n}}(Y;E,F)[[\epsilon]],$
  $n=\dim(Y),$ by $\Psi^{m}_{\suss{2n}}(Y;E,F).$
\end{definition}
\noindent The quotient here corresponds formally to setting 
\begin{equation}
\epsilon^{n+1}=0.
\label{fipomb2.40}\end{equation}

\begin{proposition}(Essentially from \cite{perdet}) The group
\begin{multline}
G^{-\infty}_{\suss{2n}}(Y;F)=\big\{\Id +S;S\in
\Psi_{\suss{2n}}^{-\infty}(Y;F), \\
\exists\ (\Id+S)^{-1}\in \Psi_{\suss{2n}}^{0}(Y;F)\big\},
\label{detb.12}\end{multline}
with composition given by the $*$-product, admits a determinant homomorphism
\begin{equation}
\det:G^{-\infty}_{\suss{2n}}(Y;F)\longrightarrow\bbC,\
\det(A*B)=\det(A)\det(B),
\label{fipomb2.37}\end{equation}
given by
\begin{equation}
  \det(B)=\exp \left( \int_{[0,1]}\gamma^{*}\alpha_{n}\right)
\label{detb.13}\end{equation}
where $\alpha_{n}$ is the coefficient of $\epsilon^{n}$ in the 1-form
$\Tr(A^{-1}*dA)$ and $\gamma:[0,1]\longrightarrow G^{-\infty}_{\suss{2n}}(Y;F)$ is 
any smooth path with $\gamma(0)=\Id$ and $\gamma(1)=B.$
\label{detb.11}\end{proposition}

\begin{proof} In \cite{perdet} 
the determinant is defined via \eqref{detb.13} for
the full formal power series algebra with $*$-product. Since the
1-form $\alpha _n$ only depends on the term of order $n$ in the formal 
power series, and this term for a product only depends on the first $n$
terms of the factors, we can work in the quotient and \eqref{fipomb2.37} 
follows.
\end{proof}

For the group $G^{-\infty}_{\suss{2}}(X;E),$ the form $\alpha_{2}$ can be
computed explicitly.

\begin{proposition}\label{detmult.19} On $G^{-\infty}_{\suss{2}}(X;E)$
\begin{multline}
\alpha_{2}=i\pi d\mu (a) -\frac{1}{4\pi i}\int_{\bbR^{2}} \Tr \big(
(a_{0}^{-1}\frac{\pa a_{0}}{\pa t})(a_{0}^{-1}\frac{\pa a_{0}}{\pa
  \tau})a_{0}^{-1}da_{0} \\ 
-(a_{0}^{-1}\frac{\pa a_{0}}{\pa \tau})(a_{0}^{-1}\frac{\pa
  a_{0}}{\pa t})a_{0}^{-1}da_{0}\big) dt d\tau ,  
\label{fipomb2.60}\end{multline}
where  
\begin{equation}
\mu(a)=\frac{1}{2\pi^{2}i}\int_{\bbR^{2}}\Tr(a_{0}^{-1}a_{1})dt d\tau.
\label{fipomb2.59}\end{equation}
\end{proposition}

\begin{proof} For $a=a_{0}+\epsilon a_{1}\in
G^{-\infty}_{\suss{2}}(X;E),$ the inverse $a^{-1}$ of $a$ with respect to
the $*$-product is
\begin{equation}
a^{-1}=a_{0}^{-1}-\epsilon(a_{0}^{-1}a_{1}a_{0}^{-1} -
\frac{i}{2}\{a_{0}^{-1},a_{0}\}a_{0}^{-1}) \, \,,
\label{tsc.4a}\end{equation} 
where $a_{0}^{-1}$ is the inverse of $a_{0}$  in $G^{-\infty}_{\susn{2}}(X;E)$ and 
\[
    \{a,b\}= D_{t}a D_{\tau}b - D_{\tau}a D_{t}b = \pa_{\tau}a\pa_{t}b-\pa_{t}a\pa_{\tau}b
\]
is the Poisson Bracket.  Hence,
\begin{equation}
\begin{aligned}
\Tr(a^{-1}*da) &= \Tr\bigg( a_{0}^{-1}da_{0}+\epsilon(-\frac{i}{2}
\{a_{0}^{-1},da_{0}\}-a_{0}^{-1}a_{1}a_{0}^{-1}da_{0} \\
  & \hspace{3cm} +\frac{i}{2}
\{a_{0}^{-1},a_{0}\}a_{0}^{-1}da_{0} + a_{0}^{-1}da_{1})\bigg) \\
&=\frac{1}{2\pi}\int_{\bbR^{2}}\bigg( \Tr( a_{0}^{-1}da_{0})
+\epsilon \Tr(-\frac{i}{2}
\{a_{0}^{-1},da_{0}\}\\
&\hspace{1cm}
-a_{0}^{-1}a_{1}a_{0}^{-1}da_{0} +\frac{i}{2}
\{a_{0}^{-1},a_{0}\}a_{0}^{-1}da_{0} + a_{0}^{-1}da_{1})\bigg)dtd\tau \, \,.
\end{aligned}
\label{tsc.5}\end{equation}
So
\begin{multline}
\alpha_{2}(a)=\frac{1}{2\pi}\int_{\bbR^{2}}\Tr\bigg(-\frac{i}{2}
\{a_{0}^{-1},da_{0}\}
-a_{0}^{-1}a_{1}a_{0}^{-1}da_{0}\\ 
+\frac{i}{2}
\{a_{0}^{-1},a_{0}\}a_{0}^{-1}da_{0} + a_{0}^{-1}da_{1}\bigg)dtd\tau \, \,.
\label{tsc.6}\end{multline}
On the right hand side of \eqref{tsc.6}, the first term vanishes since it
is the integral of the trace of a Poisson bracket. 
Indeed, integrating by parts one of the terms with respect to $t$ and $\tau,$
\begin{equation}
\begin{aligned}
\int_{\bbR^{2}}\Tr(\{a_{0}^{-1},da_{0}\})dt d\tau
&=\int_{\bbR^{2}}\Tr(D_{t}a_{0}^{-1}D_{\tau}(da_{0})-
         D_{\tau}a_{0}^{-1}D_{t}(da_{0}))dtd\tau \\
&=\int_{\bbR^{2}}\Tr(D_{t}a_{0}^{-1}D_{\tau}(da_{0})+
         D_{t}D_{\tau}a_{0}^{-1}(da_{0}))dtd\tau \\       
&=\int_{\bbR^{2}}\Tr(D_{t}a_{0}^{-1}D_{\tau}(da_{0})-
         D_{t}a_{0}^{-1}D_{\tau}(da_{0}))dtd\tau \\
&=0.
\end{aligned}
\label{tsc.7}\end{equation}
Hence,
\begin{equation}
\alpha_{2}(a)=\frac{1}{2\pi}\int_{\bbR^{2}}\Tr(
\frac{i}{2}\{a_{0}^{-1},a_{0}\}a_{0}^{-1}da_{0}
+a_{0}^{-1}da_{1}-a_{0}^{-1}a_{1}a_{0}^{-1}da_{0} )dt d\tau  .
\label{detmult.26}\end{equation}
The last two terms on the right combine to give $i\pi d\mu(a).$  Writing
out the Poisson bracket in terms of $t$ and $\tau$, gives \eqref{fipomb2.60}.
\end{proof}

\begin{proposition}\label{bour.3} Integration of the 1-form $\alpha_{2}$
gives an isomorphism
\begin{equation}
\phi:\pi_{1}\left(G^{-\infty}_{\suss{2}}(X;E)\right)\ni
f\longmapsto \frac{1}{2\pi i}\int_{\bbS^{1}}f^{*}\alpha_{2}\in \bbZ.
\label{bour.4}\end{equation}
\end{proposition}

\begin{proof}
By the previous proposition and Stokes' theorem,
\begin{equation}
\begin{aligned}
\int_{\bbS^{1}}f^{*}\alpha_{2}&=\frac{1}{12\pi i}\int_{\bbS^{3}}g^{*}(\Tr(
(a^{-1}da)^{3})) \\ 
  &= 2\pi i \int_{\bbS^{3}}g^{*}\beta_{2}^{odd},\ \forall\ f:\bbS^{1}\longrightarrow G^{-\infty}_{s(2)}(X;E),
\end{aligned}
\label{detmult.21}\end{equation}
where the map $g$ is defined by
\begin{equation}
\begin{matrix}
  g:\bbS^{3} \ni(s,\tau,t) \longmapsto f(s)(t,\tau)\in G^{-\infty}(X;E),
\end{matrix}
\label{detmult.22}\end{equation} 
and $\beta_{2}^{odd}=\frac{1}{6(2\pi i)^{2}}\Tr((a^{-1}da)^{3}).$ Our convention is that the orientation on
$\bbR^{2}$ is given by the symplectic form $\omega= d\tau\wedge dt$.  The
$3$-form $\beta_{2}^{odd}$ on $G^{-\infty}(X;E)$ is such that  
\begin{equation}
\begin{matrix}
 \lambda: \pi_{3}(G^{-\infty}(X;E)) \ni h\longmapsto
 \int_{\bbS^{3}}h^{*}\beta_{2}^{odd}\in \bbZ 
\end{matrix}  
\label{tsc.9}\end{equation}
is an isomorphism (see \cite{fipomb}) where
\begin{equation} 
G^{-\infty}(X;E)=\{ \Id+S; S\in\Psi^{-\infty}(X;E), ((\Id+S)^{-1}\in 
\Psi^{0}(X;E)\} \, \,.
\label{tsc.8}\end{equation}
Up to homotopy, 
$G^{-\infty}_{\susn{2}}(X;E)\cong [\bbS^{2},G^{-\infty}(X;E)]$, so
the map
$f\mapsto g$ is an isomorphism $\pi_{1}(G^{-\infty}_{\susn{2}}(X;E))\cong
\pi_{3}(G^{-\infty}(X;E)).$  Hence the proposition follows from
\eqref{detmult.21} and \eqref{tsc.9}.
\end{proof}

We may identify
\begin{equation}
\Psi_{\susn{2n}}^{m}(M/B;E,F)\subset \Psi_{\suss{2n}}^{m}(M/B;E,F)
\label{fipomb2.41}\end{equation}
as the subspace of elements independent of $\epsilon.$ For an elliptic
family $L$ with vanishing numerical index one can then consider in the same way
as above the (non-empty) principal $G^{-\infty}_{\suss{2n}}(Y_b;F_b)$ spaces
\begin{multline}
\mathcal{L}_{b}=\{ L_{b}+S_{b};
S_{b}\in\Psi_{\suss{2n}}^{-\infty}(Y_{b};E_{b},F_{b}) \\
\exists\ (L_{b}+S_{b})^{-1}\in \Psi_{\suss{2n}}^{-m}(Y_{b};F_{b},E_{b})\}
\label{detb.14}\end{multline}
forming a smooth infinite-dimensional bundle over $\fambas.$

\begin{definition}\label{detb.15} For an elliptic family $L\in
\Psi^{m}_{\suss{2n}}(M/B;E,F)$ with vanishing numerical index, the
determinant line bundle is given by
\begin{equation}
\Det(L)=\mathcal{L}\times_{G^{-\infty}_{\suss{2n}}(M/B;E)}\bbC
\label{detb.16}\end{equation}
where each fibre of $G^{-\infty}_{\suss{2n}}(M/B;E)$ acts on $\bbC$ via the
determinant of Proposition~\ref{detb.11}.
\end{definition}

\section{Cobordism invariance of the index}
\label{totdlb}

Suppose that the fibration \eqref{lodet.1} arises as the boundary of a
fibration where the fibre is a compact manifold with boundary:
\begin{equation}
\xymatrix{&\pa Z\ar@{-}[r]\ar@{_{(}->}[dl]_{\pa}
& \pa M\ar[dd]^{\pa \phi}\ar@{_{(}->}[dl]_{\pa}\\
Z\ar@{-}[r]& M\ar[dr]_{\phi}\\
&&B,}
\label{sus.8}\end{equation}
so $Z$ and $M$ are compact manifolds with boundary.  Let
$E$ and $F$ be  complex vector bundles over the manifold $M.$
Suspending the short exact sequence of Proposition~\ref{cusp.11} one
arrives at the short exact sequence
\begin{equation}
x\Psi^{m}_{\csn{k}}(Z;E,F)\longrightarrow\Psi^{m}_{\csn{k}}(Z;E,F)
\overset{I}\longrightarrow \Psi_{\susn{k+1}}^{m}(\pa Z;E,F),\ k\in\bbN.
\label{sus.30}\end{equation}

\begin{theorem}\label{sus.15} Let $L\in\Psi_{\susn{2n}}^{m}(\pa M /B;E,F)$ be an
  elliptic family of $2n$-suspended pseudodifferential operators and
  suppose that the fibration arises as the boundary of a fibration as in
  \eqref{sus.8} and that $L$ is the indicial family $L=I(P)$ of an
  elliptic family $P\in\Psi_{\csn{2n-1}}^{m}(M/B;E,F)$ of
  $(2n-1)$-suspended cusp pseudodifferential operators, then the index
  bundle of \eqref{detb.14} is trivial.
\end{theorem}

\begin{proof} Given $b\in B,$ we claim that $P_{b}$ can be perturbed by 
\begin{equation*}
Q_{b}\in \Psi_{\csn{2n-1}}^{-\infty}(M_{b};E_{b},F_{b})
\label{fipomb2.107}\end{equation*}
to become invertible.  Indeed, we may think of $P_{b}$ as a family of cusp
operators on $\bbR^{2n-1}.$ To this family we can associate the bundle
$\mathcal{I}_{b}$ over $\bbR^{2n-1}$ of invertible perturbations by elements in
$\Psi_{\cusp}^{-\infty}(M_{b};E_{b},F_{b})$.  This bundle is well-defined
in the sense that invertible perturbations exist for all $t\in \bbR^{2n-1}$ 
by Theorem 5.2 of \cite{fipomb}. The ellipticity of $P_{b}$ ensures that there 
exists $R>0$ such that $P_{b}(t)$ is invertible for $|t|\ge R.$  By the 
contractibility result of \cite{fipomb}, there exists an invertible section  
$P_{b}(t)+Q_{b}(t)$ of $\mathcal{I}_{b}$ such that $Q_{b}(t)=0$ for
$|t|>R$.  In particular $Q_{b}$ is an element of 
$\Psi_{\csn{2n-1}}^{-\infty}(M_{b};E_{b},F_{b})$, and so $P_{b}+Q_{b}$
is the desired invertible perturbation.

It follows that there exists $S_{b}\in \Psi_{\susn{2n}}^{-\infty}(\pa
Z_{b};E_{b},F_{b})$ such that $I(P_{b})=L_{b}$ is invertible. This could
also have been seen directly using $K$-theory and the cobordism invariance
of the index. In any case, this shows that the family $P$ gives rise to a
bundle $\mathcal{P}_{\csn{2n-1}}$ on the manifold $B$ with fibre at $b\in B$
\begin{multline}
\mathcal{ P}_{\csn{2n-1},b}=\{
P_{b}+ Q_{b};\ Q_{b}\in \Psi_{\csn{2n-1}}^{-\infty}(Z_{b};E_{b},F_{b}),\\
     (P_{b}+Q_{b})^{-1}\in \Psi_{cs(2n-1)}^{-k}(Z_{b};F_{b},E_{b})\}.
\label{sus.11}\end{multline} 
If we consider the bundle of groups
$G^{-\infty}_{\csn{2n-1}}(M/B;E)\longrightarrow B$ with fibre at $b\in B$
\begin{multline}
G^{-\infty}_{\csn{2n-1}}(Z_{b};E_{b})=
\{\Id + Q_{b};\ Q_{b}\in \Psi_{\csn{2n-1}}^{-\infty}(Z_{b};E_{b}),\\
(\Id+Q_{b})^{-1}\in \Psi_{\csn{2n-1}}^{0}(Z_{b};E_{b})\},
\label{sus.12}\end{multline}
then $\mathcal{P}_{cs(2n-1)}$ may be thought as a principal
$G^{-\infty}_{\csn{2n-1}}(M/B;E)$-bundle, where the group 
$G^{-\infty}_{\csn{2n-1}}(Z_{b};E_{b})$ acts on the right in the obvious
way. From \cite{fipomb} it follows $G^{-\infty}_{\csn{2n-1}}(Z_{b};F_{b}),$
is weakly contractible. Hence, $\mathcal{P}_{\csn{2n-1}}$ has a global
section defined over $B$, so is trivial as a principal
$G^{-\infty}_{\csn{2n-1}}(M/B;E)$-bundle. Taking the indicial family of
this global section gives a global section of the bundle $\mathcal{L}$
which is therefore trivial as a principal
$G^{-\infty}_{\suss{2n}}(M/B;E)$-bundle.
\end{proof}

As an immediate consequence, the determinant bundle of a $2n$-suspended
family which arises as the indicial family of elliptic cusp operators is
necessarily trivial. Indeed, it is an associated bundle to the index
bundle, which is trivial in that case. In the case of a twice-suspended
family we will give an explicit trivialization in terms of the extended
$\tau$ invariant of the elliptic cusp family. To do so we first need to
define the $\eta$ invariant in this context. As for the determinant of a
suspended family discussed in \cite{perdet} and in \S\ref{2n-sus} above,
the extended $\eta$ invariant is only defined on the $\star$-extended
operators which we discuss first.

\section{Suspended cusp $\star$-algebra}\label{tsc.0}

On a compact manifold $Z$ with boundary, consider, for a given boundary
defining function $x,$ the space of formal power series
\begin{equation}
\cA^m(Z;E)=\sum\limits_{j=0}^\infty\varepsilon ^jx^j\Psi_{\cs}^{m}(Z;E)
\label{fipomb2.86}\end{equation}
in which the coefficients have increasing order of vanishing at the boundary.
The exterior derivations $D_t$ (differentiation with
respect to the suspending parameter) and $D_{\log x}$ can be combined to
give an exterior derivative $D=(D_t,D_{\log x})$ valued in $\bbR^2.$  Here, the 
derivation $D_{\log x}$ is defined to be 
\begin{equation*}
    D_{\log x}A=\frac{d}{dz} x^{z} A x^{-z}  \big|_{z=0}. 
\label{fipomb2.103}\end{equation*}
for $A\in \Psi_{\cs}^{m}(Z;E)$.  It is such that (\cf \cite{fipomb} where a
different convention is used for the indicial family) 
\begin{equation*}
     I(D_{\log x}A)= 0, \quad I(\frac{1}{x} D_{\log x} A))= D_{\tau}I(A)
\label{fipomb2.104}\end{equation*}
where $\tau$ is the suspension variable for the indicial family.
Combining this with the symplectic form on $\bbR^2$ gives a star product 
\begin{equation}
\begin{gathered}
A*B=\sum\limits_{j,k,p}\varepsilon ^{j+k+p}\frac{i^{p}}{2^{p}p!}
              \omega(D_A,D_B)^{p}A_jB_k,\\
A=\sum\limits_{j}\varepsilon ^jA_j,\ B=\sum\limits_{k}\varepsilon ^kB_k.
\end{gathered}
\label{fipomb2.61}\end{equation}
Here, the differential operator $\omega(D_A,D_B)^p$ is first to be expanded
out, with $D_A$ being $D$ acting on $A$ and $D_B$ being $D$ acting on $B$
and then the product is taken in $\Psi_{\cs}^{m}(Z;E).$ Note that 
\begin{equation*}
D_{\log x}:x^p\Psi_{\cs}^{m}(Z;E)\longrightarrow x^{p+1}\Psi_{\cs}^{m-1}(Z;E)
\label{fipomb2.87}\end{equation*}
so the series in \eqref{fipomb2.61} does lie in the space \eqref{fipomb2.86}.
The same formal argument as in the usual case shows that this is an
associative product. We take the quotient by the ideal spanned by
$(\varepsilon x)^2$ and denote the resulting algebra
$\Psi_{\css}^{m}(Z;E).$ Its elements are sums $A+\varepsilon A',$ $A'\in
x\Psi_{\cs}^{m}(Z;E)$ and the product is just 
\begin{multline}
(A+\varepsilon A')*(B+\varepsilon B')=AB+
\varepsilon (AB'+A'B) \\
-\frac{i\varepsilon}{2}\left(D_{t}AD_{\log{x}}B-D_{\log{x}}AD_{t}B\right)
\mod \varepsilon^2x^2.  
\label{fipomb2.62}\end{multline}

The minus sign comes from our definition of the symplectic form \eqref{symp.0}.
The boundary asymptotic expansion \eqref{fipomb2.44}, now for suspended
operators, extends to the power series to give a map into triangular,
doubly-suspended, double power series 
\begin{equation}
I_*:\cA^m(Z;E)\longrightarrow\Psi_{\susn{2}}^m(\pa Z;E)[[\varepsilon x,x]].
\label{fipomb2.63}\end{equation}
To relate this more directly to the earlier discussion of star products on
the suspended algebras we take the quotient by the ideal generated by
$x^2$ giving a map 
\begin{equation}
\tI:\Psi_{\css}^{m}(Z;E)\longrightarrow\left\{a_{0}+xe+\varepsilon xa_{1},\
a_{0},\ e, a_1\in\Psi_{\susn{2}}^m(\pa Z;E)\right\}. 
\label{fipomb2.64}\end{equation}
The surjectivity of the indicial map shows that this map too is surjective
and so induces a product on the image.

\begin{proposition}\label{bour.2} The surjective map $\tI$ in
\eqref{fipomb2.64} is an algebra homomorphism for the product generated by 
\begin{equation}
a_0\tstar b_0=a_0b_0  
-\varepsilon x\frac{i}2\left(D_ta_0D_{\tau}b_0-D_{\tau} a_0D_tb_0\right),\
a_0,\ b_0\in\Psi^m_{\susn{2}}(\pa Z;E)
\label{fipomb2.66}\end{equation}
extending formally over the parameters $\varepsilon x,$ $x$ to the range in
\eqref{fipomb2.64}.
\end{proposition}

\begin{proof} First observe that in terms of the expansions
  \eqref{fipomb2.43} for $A\in\Psi_{\cs}^{m}(Z;E)$ and
  $B\in\Psi_{\cs}^{m'}(Z;E)$ at the boundary 
\begin{equation*}
\begin{gathered}
I_*(AB)=A_0B_0+x(A_0 F+ E B_0)+O(x^2), \\
A= A_{0} +xE +\mathcal{O}(x^{2}), \quad B=B_{0}+ xF +\mathcal{O}(x^{2}).
\end{gathered}
\label{fipomb2.67}\end{equation*}
It follows that for 
\begin{equation*}
A=A_{0}+xE+\varepsilon xA_{1},\
B=B_{0}+xF+\varepsilon xB_{1}
\label{fipomb2.69}\end{equation*}
the image of the product is
\begin{equation}
\tI(A*B)=\tI(A)\tI(B) 
-\varepsilon x\frac i2(D_t a_0D_{\tau}b_0-D_{\tau}a_0D_tb_0),
\label{fipomb2.68}\end{equation}
where $a_{0}= I(A_{0})$ and $b_{0}= I(B_{0}).$ This is precisely what is claimed.
\end{proof}

For any manifold without boundary $Y$ we will denote by $\Psi^m_{\tsuss{2}}(Y;E)$
the corresponding algebra with the product coming from \eqref{fipomb2.68}
so that \eqref{fipomb2.64} becomes the homomorphism of algebras 
\begin{equation}
\tI:\Psi_{\css}^{m}(Z;E)\longrightarrow\Psi^m_{\tsuss{2}}(\pa Z;E).
\label{fipomb2.78}\end{equation}
As the notation indicates, this algebra is closely related to
$\Psi^m_{\suss{2}}(Y;E)$ discussed in \S\ref{2n-sus}. Namely, by
identifying the parameter as $\epsilon =\varepsilon x$ the latter may be
identified with the quotient by the ideal
\begin{equation}
\xymatrix@1{
x\Psi^m_{\susn{2}}\ar[r]&
\Psi^m_{\tsuss{2}}(Y;E)\ar[r]^{\epsilon =\varepsilon x}&
\Psi^m_{\suss{2}}(Y;E).}
\label{fipomb2.71}\end{equation}
Similarly, for the invertible elements of order zero, 
\begin{equation}
\xymatrix@1{\widetilde{G}^0_{\susn{2}}(Y;E)\ar[r]&
G^0_{\tsuss{2}}(Y;E)\ar[r]^{\epsilon =\varepsilon x}&
G^0_{\suss{2}}(Y;E)}
\label{fipomb2.72}\end{equation}
is exact, where 
\begin{equation*}
\widetilde{G}^0_{\susn{2}}(Y;E)=\{ \Id+Q \in G^0_{\tsuss{2}}(Y;E) ; Q\in
x\Psi^0_{\susn{2}}(Y;E)\}.
\label{fipomb2.126}\end{equation*}

\section{Lifting the determinant from the boundary}\label{Lifting}

As a special case of \eqref{fipomb2.72} the groups of order $-\infty$
perturbations of the identity are related in the same way: 
\begin{multline}
\widetilde{G}^{-\infty}_{\susn{2}}(Y;E)\longrightarrow
G^{-\infty}_{\tsuss{2}}(Y;E)=\\
\left\{\Id+A_{0}+xE+\varepsilon xA_{1},\
A_{0},\ E,\ A_1\in\Psi_{\susn{2}}^{-\infty}(Y;E);
\Id+A_{0}\in G^0_{\susn{2}}(Y;E)\right\}\\
\longrightarrow 
G^{-\infty}_{\suss{2}}(Y;E),
\label{fipomb2.73}\end{multline}
where 
\begin{equation*}
\widetilde{G}^{-\infty}_{\susn{2}}(Y;E)=\{\Id+Q \in
G^{-\infty}_{\susn{2}}(Y;E) ; Q\in x\Psi^{-\infty}_{\susn{2}}(Y;E)\}.
\label{fipomb2.127}\end{equation*}
The determinant defined on the quotient group lifts to a homomorphism on
the larger group with essentially the same properties. In fact, it can be
defined directly as
\begin{equation}
\det(b)=\exp\left(\int_{0}^{1}\gamma^{*}\talpha_{2}\right),\ b\in
G^{-\infty}_{\tsuss{2}}(Y;E)
\label{fipomb2.74}\end{equation}
where $\tilde\alpha_{2}$ is the coefficient of $\varepsilon x$ in the
expansion of $a^{-1}\tstar da$ and $\gamma$ is a curve from the identity to
$b.$ Since the normal subgroup in \eqref{fipomb2.73} is affine, the larger
group is contractible to the smaller. Certainly the pull-back of
$\talpha_{2}$ to the subgroup is $\alpha_{2},$ with $\epsilon$ replaced by
$\varepsilon x,$ so Proposition~\ref{bour.3} holds for the larger group as
well. Indeed a minor extension of the computations in the proof of
Proposition~\ref{detmult.19} shows that at $a=(a_0+xe_2+\varepsilon
xa_1)\in G^{-\infty}_{\tsuss{2}}(Y;E)$
\begin{equation}
\begin{gathered}
\begin{aligned}
a^{-1}&=a_{0}^{-1}-x\left(a_0^{-1}ea_0^{-1}\right)
+\varepsilon x\left(a_{0}^{-1}a_{1}a_{0}^{-1} +
\frac{i}{2}\{a_{0}^{-1},a_{0}\}a_{0}^{-1}\right)
\end{aligned}\\
\Longrightarrow \talpha_2=\alpha_2
\end{gathered}
\label{fipomb2.76}\end{equation}
in terms of formula \eqref{fipomb2.60}.

Since the group  
\begin{equation}
G_{\css}^{-\infty}(Z;E)=\left\{\Id+A,\ A\in \Psi_{\css}^{-\infty}(Z;E);\ \exists\
(\Id+A)^{-1}\in\Id+\Psi_{\css}^{-\infty}(Z;E)\right\}
\label{fipomb2.80}\end{equation}
is homotopic to its principal part, and hence is contractible, the lift of
$d\log\det$ under $\tI$ must be exact; we compute an explicit formula for
the lift of the determinant.

\begin{theorem}
On $G^{-\infty}_{\css}(Z;F),$ 
\begin{equation}
\det(\tI(A))=e^{i\pi\eta _{\cusp}(A)}
\label{fipomb2.81}\end{equation}
where 
\begin{equation}
\eta_{\cusp}(A)=
\frac{1}{2\pi i}\int_{\bbR}
\bTr\left(A_{0}^{-1}\frac{\pa A_{0}}{\pa
  t} +\frac{\pa A_{0}}{\pa t}A_{0}^{-1} \right) dt + \mu(\tI(A)),
\label{fipomb2.82}\end{equation}
with $\mu$ defined in Proposition~\ref{detmult.19}.
\label{detmult.24}\end{theorem}

\begin{proof}
We proceed to compute $d\eta _{\cusp},$
\begin{multline}
d(\eta_{\cusp}-\mu)(A)=\frac{1}{2\pi i}\int_{\bbR} \bTr\big(
-(A_{0}^{-1}dA_{0})(A_{0}^{-1}\frac{\pa A_{0}}{\pa t}) +
A_{0}^{-1}\frac{\pa dA_{0}}{\pa t}\\ 
+\frac{\pa dA_{0}}{\pa t}A_{0}^{-1} - (\frac{\pa A_{0}}{\pa
  t}A_{0}^{-1})(dA_{0}A_{0}^{-1})\big) dt. 
\label{detmult.27}\end{multline}
Integrating by parts in the second and third terms gives
\begin{equation}
d(\eta_{\cusp}-\mu)(A)=\frac{1}{2\pi i}\int_{\bbR}
\bTr\left([A_{0}^{-1}\frac{\pa A_{0}}{\pa t},A_{0}^{-1}dA_{0}] -[\frac{\pa
    A_{0}}{\pa t}A_{0}^{-1},dA_{0}A_{0}^{-1}]\right)dt. 
\label{detmult.28}\end{equation}
Using the trace-defect formula, this becomes
\begin{multline}
d(\eta_{\cusp}-\mu)(A)=-\frac{1}{4 \pi^2}
\int_{\bbR^{2}}\Tr\big(a_{0}^{-1}\frac{\pa a_{0}}{\pa t}\frac{\pa}{\pa
  \tau}\left(a_{0}^{-1}da_{0} \right) \\ 
-\frac{\pa a_{0}}{\pa t}a_{0}^{-1}\frac{\pa}{\pa \tau}\left(da_{0}a_{0}^{-1}
\right) \big)dtd\tau , 
\label{detmult.29}\end{multline}
where $a_{0}=I(A_{0}).$ Expanding out the derivative with respect
to $\tau$ and simplifying
\begin{multline}
d(\eta_{\cusp}-\mu)(A)=\frac{1}{4\pi^{2}}\int_{\bbR^{2}} \Tr \big(
(a_{0}^{-1}\frac{\pa a_{0}}{\pa t})(a_{0}^{-1}\frac{\pa a_{0}}{\pa
  \tau})a_{0}^{-1}da_{0}\\ 
-(a_{0}^{-1}\frac{\pa a_{0}}{\pa \tau})(a_{0}^{-1}\frac{\pa
  a_{0}}{\pa t})a_{0}^{-1}da_{0} \big) dtd\tau , 
\label{detmult.30}\end{multline}
which shows that $i\pi d\eta_{\cusp}(A)=d\log{det(\tI(A))}.$
\end{proof}

Now, consider the subgroup
\begin{equation}
G_{\css,\tI=\Id}^{-\infty}(Z;E)\subset G_{\css}^{-\infty}(Z;E)
\label{fipomb2.83}\end{equation}
consisting of elements of the form $\Id+Q$ with
$Q\in\Psi^{-\infty}_{\cs}(Z;F)$ and $\tI(Q)=0.$ In particular  
\begin{equation}
\Id+Q\in G_{\css,\tI=\Id}^{-\infty}(Z;E)\Longrightarrow I(Q_0)=0.
\label{fipomb2.84}\end{equation}

\begin{proposition} In the commutative diagramme
\begin{equation}
\xymatrix{
G_{\css,\tI=\Id}^{-\infty}(Z;F)\ar[r]\ar[d]^{\ind}&
G_{\css}^{-\infty}(Z;F)\ar[r]^{\tI}\ar[d]^{\frac12\eta_{\cusp}}&
G_{\tsuss{2}}^{-\infty}(\pa Z;F)\ar[d]^{\det}\\
\bbZ\ar[r]&
\bbC\ar[r]^{\exp(2\pi i\cdot)}&
\bbC^*}
\label{detmult.31}\end{equation}
the top row is an even-odd classifying sequence for K-theory.
\label{detmult.32}\end{proposition}

\begin{proof} We already know the contractibility of the central group, and
the end groups are contractible to their principal parts, which are
classifying for even and odd K-theory respectively. For $A=\Id+Q\in
G_{\css,\tI=\Id}^{-\infty}(Z;E),$ $Q=Q_{0}+\epsilon Q_{1},$ 
and the condition $\tI(Q)=0$ reduces to $Q_{0},$ $Q_{1}\in
x^2\Psi^{-\infty}_{\cs}(Z;F)$ so are all of trace class. Then
\begin{equation}
\frac{1}{2}\eta_{\cusp}(A)= \frac{1}{2\pi i}\int_{\bbR} \Tr(A_{0}^{-1}
\frac{\pa A_{0}}{\pa t})dt
\label{detmult.33}\end{equation}
which is the formula for the odd index \eqref{fipomb2.26}.
\end{proof}

\section{The extended $\eta$ invariant}
\label{eta-inv}

Next we show that the cusp $\eta$-invariant defined in \eqref{fipomb2.82}
can be extended to a function on the elliptic invertible elements of
$\Psi_{\css}^{m}(Z;E).$ To do so the boundary-regularized trace $\bTr(A),$
defined on operators of order $-\dim(Z)-1,$ is replaced by a fully
regularized trace functional on $\Psi_{\cs}^{\bbZ}(Z;E)$ following the
same approach as in \cite{MR96h:58169}.

For $m\in\bbZ$ arbitrary and $A\in\Psi^{m}_{\cs}(Z;E),$
\begin{equation}
\frac{d^{p}A(t)}{dt^{p}}\in \Psi^{m-p}_{\cs}(Z;E), 
\label{detbundle.5}\end{equation}   
so the function 
\begin{equation}
h_{p}(t)= \bTr \left(\frac{d^{p}A(t)}{dt^{p}}\right)\in\CI(\bbR), 
\label{detbundle.6}\end{equation} 
is well-defined for $p>m+\dim(Z)+1.$ Since the regularization in the trace
functional is in the normal variable to the boundary, $h_{p}(t)$ has, as in the
boundaryless case, a complete asymptotic expansion as $t\longrightarrow\pm\infty$,
\begin{equation}
h_{p}(t) \sim \sum_{l\ge 0} h_{p,l}^{\pm} |t|^{m-p+\dim(Z)-l} .
\label{detbundle.7}\end{equation}
So
\begin{equation}
g_{p}(t)= \int_{-t}^{t}\int_{0}^{t_{p}} \cdots \int_{0}^{t_{1}}
h_{p}(r)drdt_{1}\ldots dt_{p} 
\label{detbundle.8}\end{equation}
also has an asymptotic expansion as $t\to\infty$,
\begin{equation}
g_{p}(t) \sim \sum_{j\ge0} g_{p,j} t^{m+1+\dim(Z)-j} + g_{p}'(t)+
g_{p}''(t)\log t, 
\label{detbundle.9}\end{equation}
where $g_{p}'(t)$ and  $g_{p}''(t)$ are polynomials of degree at most $p$.
Increasing $p$ to $p+1$ involves an additional derivative in
\eqref{detbundle.6} and an additional integral in \eqref{detbundle.8}. This
changes the integrand of the final integral in \eqref{detbundle.8} by a
polynomial so $g_{p+1}(t)-g_p(t)$ is a polynomial without constant
term. This justifies

\begin{definition}
The \emph{doubly regularized trace} is the continuous linear map 
\begin{equation}
\overline{\overline{\Tr}}:\Psi_{\cs}^{\bbZ}(Z;E)\longrightarrow \bbC
\label{fipomb2.85}\end{equation}
given by the coefficient of $t^{0}$ in the expansion \eqref{detbundle.9}.
\label{detbundle.10}\end{definition}

When $m<-1-\dim(Z),$ this reduces to the integral of the
boundary-regularized trace
\begin{equation}
\overline{\overline{\Tr}}(A)= \int_{\bbR} \bTr(A(t))dt.
\label{detbundle.11}\end{equation}
In general, the doubly regularized trace does not vanish on commutators.  
However, it does vanish on commutators where one factor vanishes to high
order at the boundary so the associated trace-defect can only involve
boundary terms.

The trace-defect formula involves a similar regularization of the trace
functional on the boundary for doubly suspended operators. So for a vector
bundle over a compact manifold without boundary consider
\begin{equation}
\int_{\bbR^{2}} \Tr(b) dtd\tau,\ b\in \Psi_{\susn{2}}^{m}(Y;E),\
m<-\dim(Y)-2.
\label{detbundle.12}\end{equation}
For general $b\in \Psi_{\susn{2}}^{\bbZ}(Y;E)$ set
\begin{equation}
\tilde{h}_{p}(t)= \int_{\bbR} \Tr\left( \frac{\pa^{p}b(t,\tau)}{\pa
    t^{p}} \right) d\tau,\ p>m +\dim(Y)+2.
\label{detbundle.13}\end{equation} 
As $t\to \pm\infty$, there is again a complete asymptotic expansion 
\begin{equation}
  \tilde{h}_{p}(t) \sim \sum_{l\ge 0} h_{p,l}^{\pm} |t|^{m+1+\dim(Y)-p-l}
\label{detbundle.14}\end{equation}
so
\begin{equation}
\tilde{g}_{p}(t)= \int_{-t}^{t}\int_{0}^{t_{p}} \cdots \int_{0}^{t_{1}}
\tilde{h}_{p}(r)drdt_{1}\ldots dt_{p} 
\label{detbundle.15}\end{equation}
has an asymptotic expansion as $t\to +\infty$
\begin{equation}
\tilde{g}_{p}(t) \sim \sum_{j\ge0} \tilde{g}_{p,j} t^{m+2 + \dim(Y)-j} +
\tilde{g}_{p}'(t)+ \tilde{g}_{p}''(t)\log{t} , 
\label{detbundle.16}\end{equation}
where $\tilde{g}_{p}'(t)$ and  $\tilde{g}_{p}''(t)$ are polynomials of
degree at most $p.$

\begin{proposition}\label{fipomb2.46} For $a\in\Psi_{\susn{2}}^{\bbZ}(Y)$
the \emph{regularized trace} $\rdsTr(a),$ defined as the coefficient of
$t^{0}$ in the expansion \eqref{detbundle.16} is a well-defined trace
functional, reducing to
\begin{equation}
\rdsTr(a)=\int_{\bbR^{2}} \Tr(a)dt d\tau,\Mwhen m<-\dim(Y)-2
\label{detbundle.18}\end{equation}
and it satisfies 
\begin{equation}
\rdsTr(\frac{\pa a}{\pa\tau})=0.
\label{fipomb2.88}\end{equation}
\end{proposition}

\begin{proof} That $\rdsTr(a)$ is well-defined follows from the discussion
  above. That it vanishes on commutators follows from the same arguments as
  in \cite{MR96h:58169}. Namely, the derivatives of a commutator,
  $\frac{d^p}{dt^p}[A,B],$ are themselves commutators and the sums of the
  orders of the operators decreases as $p$ increases. Thus, for large $p$
  and for a commutator, the function $\tilde{h}_{p}(t)$ vanishes. The
  identity \eqref{fipomb2.88} follows similarly.
\end{proof}

\begin{proposition}[Trace-defect formula] For $A,$ $B\in \Psi_{\cs}^{\bbZ}(Z),$
\begin{equation}
    \overline{\overline{\Tr}}([A,B])= \frac{1}{2\pi i} 
\rdsTr\left(I(A,\tau)\frac{\pa I(B,\tau)}{\pa\tau}\right)=-\frac{1}{2\pi i} 
\rdsTr\left(I(B,\tau)\frac{\pa I(A,\tau)}{\pa\tau}\right).
\label{tdfm.3}\end{equation}
\label{tdfm.2}\end{proposition}

\begin{proof}
For $p\in \bbN$ large enough, we can apply the trace-defect formula
\eqref{cusp.18} to get  
\begin{equation*}
\overline{\Tr}\left(\frac{\pa^p}{\pa t^p}[A,B]\right)=
\frac{1}{2\pi i}\int_{\bbR} \Tr\left(\frac{\pa^p}{\pa t^p} 
\left(I(A,\tau)\frac{\pa}{\pa \tau}I(B,\tau)
\right)\right)d\tau,
\label{fipomb2.47}\end{equation*}
from which the result follows.
\end{proof}

Using the regularized trace functional, $\mu$ may be extended from
$G^{-\infty}_{\tsuss{2}}(Y;E,F)$ to $G^{m}_{\tsuss{2}}(Y;E)$ by setting
\begin{equation}
\mu(a)= \frac{1}{2\pi^{2}i}\rdsTr(a_{0}^{-1}a_{1}),\
a=(a_0+xe+ \varepsilon x a_1)\in G^{m}_{\tsuss{2}}(Y;E).
\label{fipomb2.49}\end{equation}

\begin{proposition}\label{detbundle.20}
For $A=A_{0}+\varepsilon x A_{1}\in G_{\css}^{m}(Z;E,F),$ the set of 
invertible elements of $\Psi_{\css}^{m}(Z;E,F),$
\begin{equation}
   \etac(A) :=\frac{1}{2\pi i}\overline{\overline{\Tr}}
\left( A_{0}^{-1}\frac{\pa
     A_{0}}{\pa t} +\frac{\pa A_{0}}{\pa t}A_{0}^{-1} \right)+\mu(\tI(A)), 
\label{fipomb2.48}\end{equation}
is log-multiplicative under composition 
\begin{equation}
      \etac(A*B)= \etac(A) + \etac(B),\ \forall\ A\in
      G_{\css}^{m}(Z;E,F),\ B\in G_{\css}^{m'}(Z;F,G).
\label{fipomb2.50}\end{equation}
\end{proposition}

\begin{proof}
If $a=\tI(A,\tau)$ and $b=\tI(B,\tau)$ denote the associated boundary
operators, a straightforward calculation shows that 
\begin{equation}
\mu(a * b)=\mu(a)+\mu(b)
-\frac{1}{4\pi^{2}}\rdsTr
(b_{0}^{-1}a_{0}^{-1}\{a_{0},b_{0}\}).
\label{detbundle.21}\end{equation}
On the other hand,
\begin{multline}
\overline{\overline{\Tr}}\left(
(A_{0}B_{0})^{-1}\frac{\pa(A_{0}B_{0})}{\pa t} +\frac{\pa(A_{0}B_{0})}{\pa
  t}(A_{0}B_{0})^{-1} \right)= \\ 
\overline{\overline{\Tr}}\left( A_{0}^{-1}\frac{\pa A_{0}}{\pa
  t} +\frac{\pa A_{0}}{\pa t}A_{0}^{-1} \right)
+\overline{\overline{\Tr}}\left( B_{0}^{-1}\frac{\pa B_{0}}{\pa
  \tau} +\frac{\pa B_{0}}{\pa \tau}B_{0}^{-1} \right) +\alpha ,
\label{detbundle.22}\end{multline}
where
\begin{equation} 
\alpha =\overline{\overline{\Tr}}\left( [
  B_{0}^{-1}A_{0}^{-1}\frac{\pa A_{0}}{\pa t}, B_{0}] + [A_{0}, \frac{\pa
    B_{0}}{\pa t}B_{0}^{-1}A_{0}^{-1}] \right).
\end{equation} 
Using the trace-defect formula \eqref{tdfm.3},
\begin{equation}
\begin{aligned}
\alpha=&\frac{1}{2\pi i} \rdsTr\left(
b_{0}^{-1}a_{0}^{-1}\frac{\pa a_{0}}{\pa t}\frac{\pa b_{0}}{\pa \tau}
- \frac{\pa a_{0}}{\pa\tau}\left( \frac{\pa b_{0}}{\pa
  t}b_{0}^{-1}a_{0}^{-1}\right) \right)\\
=&-\frac{1}{2\pi i}
\rdsTr\left(b_{0}^{-1}a_{0}^{-1}\{a_{0},b_{0}\}\right).
\end{aligned} 
\label{fipomb2.89}\end{equation}
Combining \eqref{detbundle.21}, \eqref{detbundle.22} and \eqref{fipomb2.89}
gives \eqref{fipomb2.50}.
\end{proof}

\section{Trivialization of the determinant bundle}\label{Tdb.0}

In \S\ref{2n-sus} the determinant bundle is defined for a family of elliptic,
doubly-suspended, pseudodifferential operators on the fibres of a fibration
of compact manifolds without boundary. When the family arises as the
indicial family of a family of once-suspended elliptic cusp
pseudodifferential operators on the fibres of fibration, \eqref{sus.8},
the determinant bundle is necessarily trivial, following the discussion in
\S\ref{totdlb}, as a bundle associated to a trivial bundle.

\begin{theorem} If $P\in\Psi_{\cs}^m(\famrel;E,F)$ is an elliptic family of
once-suspended cusp pseudodifferential operators and $\mathcal{P}$ is the
bundle of invertible perturbations by elements of
$\Psi_{\css}^{-\infty}(\famrel;E,F)$ then the $\tau$ invariant  
\begin{equation}
\tau=\exp(i\pi\etac):\mathcal{P}\longrightarrow \bbC^*
\label{fipomb2.90}\end{equation}
descends to a non-vanishing linear function on the determinant bundle of
the indicial family
\begin{equation}
\tau:\Det(I(P))\longrightarrow \bbC.
\label{fipomb2.91}\end{equation}
\label{detbundle.26}\end{theorem}

\begin{proof} As discussed in \S\ref{totdlb}, the bundle $\mathcal{P}$ has
  non-empty fibres and is a principal 
  $G^{-\infty}_{\css}(\famrel;E)$-bundle. The cusp eta invariant, defined by
  \eqref{fipomb2.48} is a well-defined function 
\begin{equation}
\etac:\mathcal{P}\longrightarrow \bbC.
\label{fipomb2.92}\end{equation}
Moreover, under the action of the normal subgroup
  $G^{-\infty}_{\css,\tI=\Id}(\famrel;E)$ in \eqref{detmult.31} it follows
that the exponential, $\tau,$ of $\etac$ in \eqref{fipomb2.90} is
  constant. Thus  
\begin{equation}
\tau:\mathcal{P}'\longrightarrow \bbC^*
\label{fipomb2.93}\end{equation}
is well-defined where
  $\mathcal{P}'=\mathcal{P}/G^{-\infty}_{\css,\tI=\Id}(\famrel;E)$ is the 
quotient bundle with fibres which are principal spaces for the action of
  the quotient group $G^{-\infty}_{\tsuss{2}}(\pa;E)$ in
  \eqref{fipomb2.31}. In fact $\tI$ identifies the fibres of $\mathcal{P}'$ with
  the bundle $\mathcal{L}$ of invertible perturbations by
  $\Psi^{-\infty}_{\tsuss{2}}(\pa Z;E,F)$ of the indicial family of the
  original family $P,$ so
\begin{equation}
\tau:\mathcal{L}\longrightarrow\bbC^*.
\label{fipomb2.94}\end{equation}
Now, the additivity of the cusp $\eta$ invariant in \eqref{detbundle.21}
and the identification in \eqref{detmult.31} of $\tau$ with the determinant
on the structure group shows that $\tau$ transforms precisely as a linear
function on $\Det(I(P)):$
\begin{equation*}
\tau:\Det(I(P))\longrightarrow \bbC.
\label{fipomb2.95}\end{equation*}
\end{proof}

\section{Dirac families}\label{Diracfami}

In this section, we show that Theorem~\ref{detbundle.26} can be interpreted
as a generalization of a theorem of Dai and Freed in \cite{dai-freed1} for
Dirac operators defined on odd dimensional Riemannian manifolds with
boundary.  This essentially amounts to two things.  First, that the eta
functional defined in\eqref{fipomb2.48} corresponds to the usual eta
invariant in the Dirac case, which is established in
Proposition~\ref{dirac.12} below.  Since we are only defining this eta
functional for invertible operators, $e^{i\pi\eta}$ really corresponds to
the $\tau$ functional which trivializes the inverse determinant line bundle
in \cite{dai-freed1}. Secondly, that in the Dirac case, the determinant
bundle $\det(I(P))$ is isomorphic to the determinant line bundle of the
associated family of boundary Dirac operators, which is the content of
Proposition~\ref{28.07.05.1} below.

As a first step, let us recall the usual definition of the eta function 
on a manifold with boundary (see \cite{Mueller4}).
Let $X$ be a Riemannian manifold with nonempty boundary $\pa  X=Y$.  Near the 
boundary, suppose that the Riemannian metric is of product type, so there is
a neighborhood $Y\times [0,1)\subset X$ of the boundary in which the metric 
takes the form
\begin{equation}
g= du^{2} + h_{Y}
\label{dirac.1}\end{equation}
where $u\in [0,1)$ is the coordinate normal to the boundary and $h_{Y}$ is the 
pull-back of a metric on $Y$ via the projection $Y\times[0,1)\longrightarrow
Y.$  Let $S$ be a Hermitian vector bundle over $X$ and let
$D:\mathcal{C}^{\infty}(X,S)\longrightarrow \mathcal{C}^{\infty}(X,S)$ be a
first order elliptic differential operator on $X$ which is formally
selfadjoint with respect to the inner product defined by the fibre metric
of $S$ and the metric on $X.$ In the neighborhood $Y\times [0,1)\subset X$
of the boundary described above, assume that the  operator $D$ takes the form 
\begin{equation}
  D=\gamma \left( \frac{\pa}{\pa u} + A \right)
\label{dirac.2}\end{equation}
where $\gamma:S\big|_{Y}\longrightarrow S\big|_{Y}$ is a
bundle isomorphism and $A:\mathcal{C}^{\infty}(Y,S\big|_{Y})\longrightarrow  
\mathcal{C}^{\infty}(Y,S\big|_{Y})$ is a first order elliptic operator on
$Y$ such that
\begin{equation}
  \gamma^{2}=-\Id,\ \gamma^{*}=-\gamma,\ A\gamma= -\gamma A,\ A^{*}=A.
\label{dirac.3}\end{equation}   
Here, $A^{*}$ is the formal adjoint of $A.$ Notice in particular that this 
includes the case of a compatible Dirac operator when $S$ is a Clifford module 
and $\gamma= cl(du)$ is the Clifford multiplication by $du.$ If
$\ker{A}=\{0\}$, consider the spectral boundary condition
\begin{equation}
\varphi \in \mathcal{C}^{\infty}(X,S),\ \Pi_{-}(\varphi\big|_{Y})=0,
\label{dirac.4}\end{equation}
where $\Pi_{-}$ is the projection onto the positive spectrum of $A.$ In the 
case where $\ker{A}\ne \{0\},$ a unitary involution $\sigma:\ker{A}
\longrightarrow \ker{A}$ should be chosen such that
$\sigma\gamma=-\gamma\sigma$ (such an involution exists), and the
boundary condition is then modified to
\begin{equation}
\varphi \in \mathcal{C}^{\infty}(X,S),\ (\Pi_{-}+P_{-})(\varphi\big|_{Y})=0,
\label{dirac.5}\end{equation}
where $P_{-}$ is the orthogonal projection onto $\ker{(\sigma + \Id)}.$ The 
associated operator $D_{\sigma}$ is selfadjoint and has pure point spectrum.  
For this operator, the eta invariant is
\begin{equation}
\eta_{X}(D_{\sigma})=\frac{1}{\sqrt{\pi}}\int_{0}^{\infty} s^{-\frac{1}{2}}
\Tr(D_{\sigma}e^{-sD_{\sigma}^{2}})ds.
\label{dirac.6}\end{equation}

To make a link with the cusp calculus, we need to enlarge $X$ by attaching
the half-cylinder $\bbR^{+}\times Y$ to the boundary $Y$ of $X.$ The
product metric near the boundary extends to this half-cylinder, which makes
the resulting manifold a complete Riemannian manifold. Similarly, the
operator $D$ has a natural extension to $M$ using its product structure
near the boundary. Denote its $L^{2}$ extension (on $M$) by $\mathcal{D}.$
The operator $\mathcal{D}$ is selfadjoint. The eta invariant of
$\mathcal{D}$ is
\begin{equation}
\eta_{M}(\mathcal{D})= \frac{1}{\sqrt{\pi}}\int_{0}^{\infty} s^{-\frac{1}{2}}
\int_{M} tr(E(z,z,s))dz ds,
\label{dirac.7}\end{equation} 
where $E(z_{1},z_{2},s)$ is the kernel of
$\mathcal{D}e^{-s\mathcal{D}^{2}}.$ One of the main result of
\cite{Mueller4} is to establish a correspondence between the eta invariants
\eqref{dirac.6} and \eqref{dirac.7}.

\begin{theorem}[\textbf{M\"{u}ller}]\label{dirac.8}
Let $D:\mathcal{C}^{\infty}(X,S)\longrightarrow \mathcal{C}^{\infty}(X,S)$ be a 
compatible Dirac operator which, on a neighborhood $Y\times[0,1)$ of $Y$ in
$X,$ takes the form \eqref{dirac.2}, let $C(\lambda):\ker{A}\longrightarrow
\ker{A}$ be the associated scattering matrix (see \cite{Mueller4} for a
definition) in the range $|\lambda|<\mu_{1}$, where $\mu_{1}$ is the
smallest positive eigenvalue of $A$ and put $\sigma=C(0),$ then 
\begin{equation*}
\eta_{X}(D_{\sigma})=\eta_{M}(\mathcal{D}).
\label{fipomb2.105}\end{equation*}
\end{theorem}

Now, on $M,$ it is possible to relate $D$ to a cusp operator. Extending the
variable $u$ to the negative reals gives a neighborhood
$Y\times(-\infty,1)\subset M$ of $\pa X$ in $M.$ The variable
\begin{equation}
  x=-\frac{1}{u} 
\label{dirac.9}\end{equation}
takes value in $(0,1)$ and by extending it to $x=0,$ gives a manifold with
boundary $\overline{M},$ with $x$ as a boundary defining function so fixing
a cusp structure. Denote by $D_{c}$ the natural extension of $D$ to
$\overline{M}.$ Near the boundary of $\overline{M},$
\begin{equation}
D_{c}=\gamma\left( x^{2}\frac{\pa}{\pa x}+A\right) 
\label{dirac.10}\end{equation}
and so is clearly a cusp differential operator. If $S=S^{+}\oplus S^{-}$ is
the decomposition of $S$ as a superspace, then
\begin{equation}
  \hat{D}_{\cs}(t)= D_{c} +it \in \Psi^{1}_{\css}(\bar{M};S)
\label{dirac.11}\end{equation}
is a suspended cusp operator, where there are no $x$ and $\varepsilon x$
terms. When $D_{c}$ is invertible, $\hat{D}_{\cs}$ is invertible as well and 
$\eta_{\cusp}(\hat{D}_{\cs})$ is well-defined.

\begin{proposition}\label{dirac.12}
Let $X,Y,M,\overline{M}$ be as above and let $D$ be a compatible Dirac
operator for some Clifford module $S$ on $X,$ which, on a neighborhood 
$Y\times[0,1)$ of $Y$ takes the form \eqref{dirac.2}, suppose that $D$ is
invertible, and let $D_{c}$ be its extension to $\overline{M},$ then 
\begin{equation*}
\eta_{X}(D_{\sigma})=\etac(\hat{D}_{\cs}),
\label{fipomb2.106}\end{equation*}
where $\hat{D}_{\cs}=D_{c}+it \in \Psi^{1}_{\css}(\overline 
{M};S)$ and $\sigma$ is trivial since $A$ is invertible.
\end{proposition}

\begin{proof} By the theorem of M\"{u}ller, it suffices to show that
$\etac(\hat{D}_{\cs})=\eta_{M}(\mathcal{D})$.  In order to do this, we
closely follow the proof of Proposition~5 in \cite{MR96h:58169}, which is
the same statement but in the case of a manifold without boundary.

Let $E(z_{1},z_{2},s)$ denote the kernel of
$\mathcal{D}e^{-s\mathcal{D}^{2}}$, where $\mathcal{D}$ is the $L^{2}$
extension of $\eth$ on $M$.  In \cite{Mueller4}, it is shown  that $\tr(
E(z,z,s))$ is absolutely integrable on $M$, so set
\begin{equation}
h(s)=\int_{M} \tr (E(z,z,s)) dz \, , s\in [0,\infty).
\label{dirac.13}\end{equation}
Then, (see \cite{Mueller4}) for $n=\dim (X)$ even,
$h(s)\in\mathcal{C}^{\infty}([0,\infty)),$ while for $n=\dim (X)$ odd, 
$h(s)\in s^{\frac{1}{2}}C^{\infty}([0,\infty))$.  Moreover, since
 $\ker\mathcal{D}=\{0\}$, $h$ is exponentially decreasing as $s\to
+\infty.$ As in \cite{MR96h:58169}, consider
\begin{equation}
g(v,t)=\int_{v}^{\infty}e^{-st^{2}}h(s)ds \, ,v\ge 0.
\label{dirac.14}\end{equation}  
This is a smooth function of $v^{\frac{1}{2}}$ in $v\ge 0$ and $t\in \bbR$, 
and as $|t|\to\infty$, it is rapidly decreasing if $v>0$.  From the fact 
that $h(s)\in C^{\infty}([0,\infty))$ for $n$ even, 
$h(s)\in s^{\frac{1}{2}}C^{\infty}([0,\infty))$ for $n$ odd, and the
exponential decrease, we get
\begin{equation}
 \left| \left(t \frac{\pa}{\pa t}\right)^{p}g(v,t)\right|\le 
\frac{C_{p}}{1+t^{2}} \, \,\, \, , v\ge 0, t \in \bbR.
\label{dirac.15}\end{equation}
So $g$ is uniformly a symbol of order $-2$ in $t$ as $v$ approaches $0.$
In fact, when $n$ is odd, it is uniformly a symbol of order $-3.$
Now,using the identity
\begin{equation}
 1= \frac{1}{\sqrt{\pi}}\int_{-\infty}^{+\infty}s^{\frac{1}{2}}\exp(-st^{2})
dt \, \, \, , s>0,
\label{dirac.16}\end{equation}
$\eta_{M}(\mathcal{D})$ may be written as a double integral
\begin{equation}
\eta_{M}(\mathcal{D})=\frac{1}{\sqrt{\pi}}
\int_{0}^{\infty}s^{-\frac{1}{2}}h(s)ds = 
 \frac{1}{\pi} \int_{0}^{\infty}\left( \int_{-\infty}^{+\infty}e^{-st^{2}}
h(s)dt \right) ds .
\label{dirac.17}\end{equation}
The uniform estimate \eqref{dirac.15} allows the limit and integral to be
exchanged so
\begin{equation}
\eta_{M}(\mathcal{D})= \frac{1}{\pi} \lim_{v\to 0}\int_{-\infty}^{+\infty} g(v,t)
dt.
\label{dirac.18}\end{equation}
For $p\in \bbN_{0}$
\begin{equation}
 g_{p}(v,t):= \int_{-t}^{t}\int_{0}^{t_{p}}\cdots 
\int_{0}^{t_{1}}  \frac{\pa^{p}}{\pa r^{p}}g(v,r)dr dt_{1}\ldots 
dt_{p} ,
\label{dirac.19}\end{equation}
has a uniform asymptotic expansion as $t\to\infty$ and $\eta_{M}(\mathcal{D})$ is
just the limit as $v\to 0$ of the coefficient of $t^{0}$ in this expansion.  

The kernel $E(z_{1},z_{2},s)$ can also be thought as the kernel of
$D_{c}e^{-s D_{c}},$ a cusp operator of order $-\infty$ on $\overline{M}.$ This
can be checked directly from the explicit construction of
$E(z_{1},z_{2},s)$ given in \cite{Mueller4}. Therefore,
\begin{equation}
h(s)=\bTr(D_{c}e^{-sD_{c}^{2}}) , 
\label{dirac.20}\end{equation}
where $\bTr$ is the regularized trace defined in \cite{fipomb}.  Note however 
that in this case, it is just the usual trace, that is, the integral of the 
kernel along the diagonal, since the residue trace vanishes.  Consider now the 
(cusp product-suspended) operator
\begin{equation}
\hat{A}(t)=\int_{0}^{\infty}e^{-st^{2}}D_{c} e^{-s D_{c}^{2}}ds= 
\frac{D_{c}}{t^{2}+D_{c}^{2}}  \, \,.
\label{dirac.21}\end{equation}

Then, $\overline{\overline{\Tr}}(\hat{A})$ is the coefficient of $t^{0}$ in the 
asymptotic expansion as $t\to \infty$ of
\begin{equation}
\int_{-t}^{t}\int_{0}^{t_{p}}\cdots \int_{0}^{t_{1}}  
\bTr(\frac{d^{p}}{d r^{p}}\hat{A}(r))dr dt_{1}\ldots dt_{p}
\label{dirac.22}\end{equation}
for $p>n=dim(\overline{M})$ and
\begin{equation}
\begin{aligned}
\bTr\left( \frac{d^{p}}{dt^{p}}\hat{A}(t)\right) &= \bTr 
\left(\frac{d^{p}}{dt^{p}}\int_{0}^{\infty}e^{-st^{2}}
D_{c} e^{-s D_{c}^{2}} ds\right) \\
     &= \frac{d^{p}}{dt^{p}}\int_{0}^{\infty}e^{-st^{2}}
\bTr(D_{c} e^{-s D_{c}^{2}}) ds \\
     &=\frac{d^{p}}{dt^{p}}g(0,t),
\end{aligned}
\label{dirac.23}\end{equation}
so $\eta_{M}(\mathcal{D})=\frac{1}{\pi}
\overline{\overline{\Tr}}(\hat{A}).$ Instead of $\hat{A},$ consider
\begin{equation}
\hat{B}(t)=\int_{0}^{\infty}e^{-st^{2}}(D_{c}-it)e^{-s D_{c}^{2}}ds=
\frac{1}{it+ D_{c}}=(\hat{D}_{\cs})^{-1} \, \, .
\label{dirac.24}\end{equation}
Since $\hat{B}(t)-\hat{A}(t)$ is odd in $t$, 
$\overline{\overline{\Tr}}(\hat{B}-\hat{A})=0$, so finally
\begin{equation}
\begin{aligned}
\eta_{M}(\mathcal{D}) &=\frac{1}{\pi}\overline{\overline{\Tr}}((
\hat{D}_{\cs})^{-1})=
\frac{1}{2\pi i}\overline{\overline{\Tr}}\left( \frac{\pa}{\pa t}
(\hat{D}_{\cs})(\hat{D}_{\cs})^{-1} +(\hat{D}_{\cs})^{-1}
\frac{\pa}{\pa t}(\hat{D}_{\cs})
\right) \\
     &=\etac(\hat{D}_{\cs})  .
\end{aligned}
\label{dirac.25}\end{equation}
\end{proof}

Let $\eth$ be some compatible Dirac operator as in Proposition~\ref{dirac.12}.
Then near the boundary of $\overline{M}$, 
its cusp version $\eth_{c}$ takes the form
\begin{equation}
\eth_{c}=\gamma(x^{2}\frac{\pa}{\pa x} +A).
\label{detdirac.19}\end{equation}
Here, it is tacitly assumed that near
the boundary, $S$ is identified with the pull-back of 
$S\big|_{\pa\overline{M}}$ via the projection 
$\pa\overline{M}\times[0,1)\longrightarrow\pa\overline{M}$.  Since the map
\begin{equation}
T^{*}(\pa\overline{M})\ni \xi \longmapsto cl(du)cl(\xi)\in \Cl(\overline{M}),
\ \gamma=cl(du), \quad x=-\frac{1}{u},
\label{detdirac.20}\end{equation}
extends to an isomorphism of algebras
\begin{equation}
\Cl(\pa\overline{M}) \longrightarrow\Cl^{+}(\overline{M})\big|_{\pa\overline{M}},
\label{detdirac.21}\end{equation}
where $\Cl(\pa\overline{M})$ and $\Cl(\overline{M})$ are the Clifford algebras
of $\pa\overline{M}$ and $\overline{M}$, this gives an action of $\Cl(\pa
\overline{M})$ on $S^{0}=\left.S^{+}\right|_{\pa\overline{M}}$.  If 
$\nu^{+}:S^{+}_{\pa\overline{M}}\longrightarrow
S^{0}$ denotes this identification, $S^{-}_{\pa
\overline{M}}$ can be identified with $S^{0}$ via the map
\begin{equation}
\nu^{-}=\nu^{+}\circ cl(du):\left.S^{-}\right|_{\pa\overline{M}}
\longrightarrow S^{0}.
\label{detdirac.22}\end{equation}
The combined identification $\nu:\left.S\right|_{\pa\overline{M}}
\longrightarrow S^{0}\oplus S^{0}$ allows us to write $\eth_{c}$ and 
$\gamma$ as
\begin{equation}
\eth_{c}= \begin{pmatrix}
               0 & \eth_{0}+x^{2}\frac{\pa}{\pa x} \\
               \eth_{0}-x^{2}\frac{\pa}{\pa x} & 0   
                 \end{pmatrix}\
\gamma = \begin{pmatrix}
               0 & 1 \\
               -1 & 0   
                 \end{pmatrix}
\label{detdirac.23}\end{equation}
acting on $S^{0}\oplus S^{0},$ where $\eth_{0}$ is the 
Dirac operator associated to the $\Cl(\pa\overline{M})$-module $S^{0}.$ If
instead we decompose the bundle $S^{0}\oplus S^{0}$ in terms of the $\pm i$
eigenspaces $S^{\pm}$ of $\gamma,$ then $\eth_{c}$ and $\gamma$ take the form
\begin{equation}
\eth_{c}= \begin{pmatrix}
               ix^{2}\frac{\pa}{\pa x} & \eth^{-}_{0} \\
               \eth^{+}_{0} & -ix^{2}\frac{\pa}{\pa x}   
                 \end{pmatrix},\ 
\gamma = \begin{pmatrix}
               i & 0 \\
               0 & -i   
                 \end{pmatrix}
\label{detdirac.33}\end{equation} 
with $\eth^{\pm}_{0}=\pm i\eth_{0}$, so that $\eth^{+}_{0}$ and $\eth^{-}_{0}$ 
are the adjoint of each other. In this notation, the suspended operator 
$\hat{\eth}_{\cs}(t)$ can be written as 
 \begin{equation}
\hat{\eth}_{\cs}(t)= \begin{pmatrix}
                      it + ix^{2}\frac{\pa}{\pa x} & \eth^{-}_{0} \\
                       \eth^{+}_{0} & it -ix^{2}\frac{\pa}{\pa x}   
                       \end{pmatrix}.
\label{detdirac.34}\end{equation}
Thus, its indicial operator $\hat{\eth}_{s(2)}(t,\tau)$ is
\begin{equation}
\hat{\eth}_{s(2)}(t,\tau)=\left. e^{\frac{i\tau}{x}}\hat{D}_{c}(t)
e^{-\frac{i\tau}{x}}\right|_{x=0}=\begin{pmatrix}
                      it - \tau & \eth^{-}_{0} \\
                       \eth^{+}_{0} & it +\tau   
                       \end{pmatrix}.
\label{detdirac.35}\end{equation}
Note that 
\begin{equation}
\hat{\eth}_{s(2)}^{*}\hat{\eth}_{s(2)}=
\begin{pmatrix}
                       t^{2}+\tau^{2}+\eth^{-}_{0}\eth^{+}_{0} & 0 \\
                         0       & t^{2}+\tau^{2}+\eth^{+}_{0}\eth^{-}_{0}   
                       \end{pmatrix}
\label{detdirac.36}\end{equation}
which is invertible everywhere except possibly at $t=\tau=0$ so
$\hat{\eth}_{s(2)}$ is invertible for every $t,\tau \in\bbR$ if and only if
$\eth^{+}_{0}$ (and consequently $\eth^{-}_{0}$) is invertible. 

Now, we wish to relate the determinant bundle associated to the family
$\hat{\eth}_{s(2)}$ with the determinant bundle of the boundary Dirac
family $\eth_{0}^{+}$ using the periodicity of the determinant line bundle
discussed in \cite{perdet}. For Dirac operators on a closed manifold, this
can be formulated as follows.

\begin{proposition} If $\eth_{0}^{+}\in \operatorname{Diff}^{1}(N/B; S)$
is a family of Dirac type operators on a fibration of closed manifold
$N\longrightarrow B$ with vanishing numerical index then the determinant
bundle of $\eth^{+}_{0}$ is naturally isomorphic to the determinant bundle
of the associated family of twice suspended operators 
\begin{equation*}
\eth_{s(2)}(t,\tau)=\begin{pmatrix}
it-\tau & \eth_{0}^{-} \\ 
\eth_{0}^{+} & it+\tau
\end{pmatrix}
\in \Psi_{\susn{2}}^{1}(N/B;S\oplus S)
\label{fipomb2.55}\end{equation*}
where $\eth_{0}^{-}=(\eth_{0}^{+})^{*}.$ 
\label{28.07.05.1}\end{proposition}

In \cite{perdet}, this periodicity is formulated in terms of product-suspended
operators instead of suspended operators, since in general,
given $P\in \Psi^{1}(Y;S)$, the family of operators 
\begin{equation*}
P_{s(2)}(t,\tau)=\begin{pmatrix}
it-\tau & P^{*} \\ 
P & it+\tau
\end{pmatrix}
\label{fipomb2.56}\end{equation*}
is a twice product-suspended operator but not a 
suspended operator unless $P$ is a differential operator. In this latter
case, which includes Dirac operators, the periodicity of the determinant
line bundle can be formulated using only suspended operators.

\section{Generalization to product-suspended operators}\label{asg.0}

As written, Theorem~\ref{detbundle.26} applies to elliptic families of
once-suspended cusp pseudodifferential operators. As we discussed in
\S\ref{Diracfami}, this includes the result of Dai and Freed
\cite{dai-freed1} for a family of self-adjoint Dirac operators $D$ on a
manifold with boundary by passing to the associated family of elliptic cusp
operators and then to the elliptic family of once-suspended cusp operators
$D+it,$ the suspension parameter being $t.$

More generally, one can consider the case of an arbitrary elliptic family
of first order self-adjoint cusp pseudodifferential operators $P.$ Then
$P+it$ is not in general a once-suspended family of cusp operators. Instead
we pass to the larger algebra, and related modules $\Psi_{\cps}^{k,l}(M/B;E,F),$
of product-suspended cusp pseudodifferential operators since then 
\begin{equation}
P+it\in\Psi_{\cps}^{1,1}(M/B;E).
\label{fipomb2.108}\end{equation}
The ellipticity of $P$ again implies the corresponding full ellipticity of
$P+it.$

Enough of the properties of suspended (cusp) operators extend to the
product-suspended case to allow the various definitions of regularized
traces and the eta invariant to carry over to the more general case.  In
this context, the proof of Proposition~\ref{detbundle.20} still applies, so
the eta invariant is also multiplicative under composition of invertible
fully elliptic cusp product-suspended operators.

So, given an elliptic family of cusp product-suspended operators (see
Definition~\ref{elliptic.1} in the Appendix), consider the bundle
$\mathcal{P}$ of invertible perturbations by elements in 
\begin{equation*}
\Psi_{\cps}^{-\infty,-\infty}(M/B;E,F)= \Psi_{\cs}^{-\infty}(M/B;E,F).
\label{fipomb2.109}\end{equation*}
Again the fibres are non-empty since the full ellipticity of $P_{b}$
implies that it is invertible for large values of the suspension parameter
$t.$ Then the same argument as in the proof of Theorem 7.1 of \cite{fipomb}
applies, using the contractibility of $G^{-\infty}_{\cusp}(M_{b};E),$ to
show the existence of an invertible perturbation $P_{b}+Q_{b}.$  At the
same time, this shows the existence of an invertible perturbation of the
indicial family $I(P_{b})\in \Psi_{\psusn{2}}^{k,l}(M_{b};E_{b},F_{b}),$
and so the associated index bundle and determinant bundle of the indicial
family are also well-defined, in the latter case using the $*$-product as
before.

Consequently, we can formulate the following generalization of
Theorem~\ref{detbundle.26} with the proof essentially unchanged.

\begin{theorem}\label{05.09.2005} If $P\in \Psi_{\cps}^{k,l}(M/B;E,F)$ is
a fully elliptic family of cusp product-suspended pseudodifferential operators
and $\mathcal{P}$ is the bundle of invertible perturbations by elements of
$\Psi_{\css}^{-\infty}(M/B:E,F)$, then the $\tau$ invariant 
\begin{equation*}
     \tau= \exp( i\pi \eta_{\cusp}):\mathcal{P}\longrightarrow \bbC^{*}
\label{fipomb2.110}\end{equation*}
descends to a non-vanishing linear function on the determinant line bundle
of the indicial family $\tau:\Det(I(P))\longrightarrow \bbC.$ 
\end{theorem}  

As a special case, Theorem~\ref{05.09.2005} includes elliptic families of
self-adjoint first order cusp pseudodifferential operators $P\in
\Psi^{1}_{\cusp}(M/B;E,F)$ by considering the cusp product suspended family 
\begin{equation*}
           P+it\in \Psi_{\cps}^{1,1}(M/B;E,F).
\label{fipomb2.111}\end{equation*}

\appendix

\section*{Appendix. Product-Suspended Operators}\label{PSO.0}

In this appendix, we will briefly review the main properties of
product-suspended pseudodifferential operators and then discuss the steps
needed to extend this notion to the case of the cusp algebra of
pseudodifferential operators on a compact manifold with boundary as used in
\S\ref{asg.0}. For a more detailed discussion on product-suspended
operators see \cite{perdet}.

For the case of a compact manifold without boundary, product-suspended
operators are, formally, generalizations of the suspended operators
by relaxing the conditions on the the full symbols. This is achieved by
replacing the radial compactification $\com{\bbR^{p}\times T^{*}X}$ by the
following blown-up version of it 
\begin{equation}
{}^{X}\overline{\bbR^{p}\times T^{*}X}=
[\overline{\bbR^{p}\times T^{*}X};\pa( \overline{\bbR^{p}}\times X)]
\label{fipomb2.112}\end{equation}
where $X$ is understood as the zero section of $T^{*}X.$ If $\rho_{r}$ and
$\rho_{s}$ denote boundary defining functions for the `old' boundary and the
`new' boundary (arising from the blow-up) then set
\begin{equation}
   S^{z,z'}({}^{X}\overline{\bbR^{p}\times T^{*}X};\hom(E,F))
=\rho_{r}^{-z}\rho_{s}^{-z'}
\CI({}^{X}\overline{\bbR^{p}\times T^{*}X};\hom(E,F)).
\label{ps.1}\end{equation}
This is the space of `full symbols' of product-suspended pseudodifferential
operators (with possibly complex multiorders). After choosing appropriate
metrics and connections, Weyl quantization gives families of operators on
$X$ which we interpret as elements of $\Psi^{z,z'}_{\psusn{p}}(X;E,F)$,
the space of product $p$-suspended operators of order $(z,z').$ However
this map is not surjective modulo rapidly decaying smoothing operators as
is the case for ordinary (suspended) pseudodifferential operators. Rather we
need to allow as a subspace 
\begin{equation}
\Psi^{-\infty,z'}_{\psusn{p}}(X;E,F)=\rho_{s}^{-z'}\CI(\com{\bbR^p}\times
X^2;\Hom(E,F),\Omega_R))\subset \Psi^{-\infty,z'}_{\psusn{p}}(X;E,F)
\label{fipomb2.113}\end{equation}
considered as smoothing operators on $X$ with parameters in $\bbR^p.$ The
image of the symbol space \eqref{ps.1} under Weyl quantization is, modulo
such terms, independent of the choices made in its definition.
The product-suspended operators form a bigraded algebra.
\begin{property}\label{prop.prod} For all $k,$ $k',$ $l$ and $l'$ and
  bundles $E,$ $F$ and $G$,
\begin{equation*}
\Psi^{k,k'}_{\psusn{p}}(X;F,G)
\circ\Psi^{l,l'}_{\psusn{p}}(X;E,F)\subset
\Psi^{k+l,k'+l'}_{\psusn{p}}(X;E,F). 
\label{fipomb2.114}\end{equation*}
\end{property}

There are two symbol maps. The usual symbol coming from the leading part of
the full symbol in \eqref{ps.1} at the `old' boundary $(B_\sigma)$ and the
`base family' which involves both the leading term of this symbol at the
`new' boundary and the leading term of the smoothing part in
\eqref{fipomb2.113}. These symbols are related by a compatibility condition
just corresponding to the leading part of the full symbol at the corner.

\begin{property}\label{prop.symb} The two symbols give short exact
sequences 
\begin{equation}
\Psi^{k-1,k'}_{\psusn{p}}(X;E,F)\longrightarrow
\Psi^{k,k'}_{\psusn{p}}(X;E,F)\overset{\sigma}\longrightarrow
\rho ^{-k'}\CI(B_\sigma ;\hom(E,F)),
\label{fipomb2.115}\end{equation}
and  
\begin{equation}
\Psi^{k,k'-1}_{\psusn{p}}(X;E,F)\longrightarrow
\Psi^{k,k'}_{\psusn{p}}(X;E,F)\overset{\beta}\longrightarrow
\Psi^{k}((X\times\bbS^{p-1})/\bbS^{p-1};E,F)
\label{fipomb2.116}\end{equation}
and the joint range is limited only by the condition 
\begin{equation}
\sigma (\beta)=\sigma |_{\pa B_{\sigma }}.
\label{fipomb2.117}\end{equation}
\end{property}

Ellipticity of $A$ is the condition of invertibility of $\sigma (A)$ and full
ellipticity is in addition the invertibility of $\beta(A).$
\begin{property}
If a fully elliptic product-suspended operator $Q\in
\Psi^{k,k'}_{\psusn{p}}(X;E,F)$ is invertible (\ie is bijective from
$\CI(X;E)$ to $\CI(X;F))$ then its inverse is an element of
$\Psi^{-k,-k'}_{\psusn{p}}(X;F,E).$
\label{ps.4}\end{property}

In \S\ref{asg.0}, we also make use of the following important
properties.
\begin{property}\label{ps.5} For all $k\in\bbZ,$ 
\begin{gather*}
\Psi_{\susn{p}}^{k}(X;E,F)\subset \Psi_{\psusn{p}}^{k,k}(X;E,F), \\
\Psi_{\susn{p}}^{-\infty}(X;E,F)= \Psi_{\psusn{p}}^{-\infty,-\infty}(X;E,F).
\end{gather*}
\end{property}

\begin{property}
If $P\in \Psi^{1}(X;E,F)$ is any first order pseudodifferential operator, then
\begin{equation*}
               P+it\in \Psi_{\psusn{1}}^{1,1}(X;E,F) 
\label{fipomb2.119}\end{equation*}
where $t$ is the suspension parameter.
\label{ps.6}\end{property}

\begin{property}  Given $Q\in \Psi_{\psusn{p}}^{k,l}(X;E,F)$,
\begin{equation*}
\left(\frac{\pa Q}{\pa t_{i}}\right)\in \Psi_{\psusn{p}}^{k-1,l-1}(X;E,F), 
\label{fipomb2.118}\end{equation*}
where $t=(t_{1},\ldots,t_{p})$ is the suspension parameter and 
$i\in \{1,\cdots,p\}.$
\label{ps.7}\end{property}

Next we extend this to a construction of cusp product-suspended 
pseudodifferential operators on a compact manifold with boundary $Z.$ Again, for
suspended cusp operators, there is a Weyl quantization map from the
appropriate space of classical symbols $\rho ^{-z}\CI(\com{\bbR^{p}\times
  {}^{\cusp}T^{*}Z};\hom(E,F))$ which is surjective modulo cusp operators
of order $-\infty$ decaying rapidly in the parameters. To capture the
product-suspended case consider the spaces of symbols analogous to \eqref{ps.1}
\begin{equation*}
\rho_{r}^{-z}\rho_{s}^{-z'}
\CI([\com{\bbR^{p}\times{}^{\cusp}T^{*}Z};\pa(\com{\bbR^p\times Z})];\hom(E,F))
\label{fipomb2.121}\end{equation*}
with the corresponding `old' and `new' boundaries. Then an element of the
space 
\begin{equation*}
A\in\Psi_{\cpsn{p}}^{z,z'}(Z;E,F)
\label{fipomb2.122}\end{equation*}
is the sum of the Weyl quantization (for the cusp algebra) of an element of
\eqref{fipomb2.121} plus an element of the residual space 
\begin{equation}
\Psi_{\cpsn{p}}^{-\infty,z'}(Z;E,F)=\rho
^{-z'}\CI(\com{\bbR^p};\Psi^{-\infty}_{\cu}(Z;E,F)).
\label{fipomb2.123}\end{equation}

Now, with this definition the properties above carry over to the boundary
setting. Property~\ref{prop.prod} is essentially unchanged. The same two
homomorphisms are defined, the symbol and base family, with the latter
taking values in families of cusp operators. In addition the indicial family
for cusp operators leads to a third homomorphism giving a short
exact sequence 
\begin{equation*}
x\Psi_{\cpsn{p}}^{m,m'}(Z;E,F)\longrightarrow \Psi_{\cpsn{p}}^{m,m'}(Z;E,F) 
     \overset{I_{\cusp}}\longrightarrow \Psi_{\psusn{p},\sus}^{k,l}(\pa Z;E,F)
\label{fipomb2.124}\end{equation*}
where the image space has the same $p$ product-suspended variables but
taking values in the suspend operators on $\pa Z.$ Since the suspended
algebra may be realized in terms of ordinary pseudodifferential operators
on $\bbR\times\pa Z$ there is no difficulty in considering these
`mixed-suspended' operators.

\begin{definition}\label{elliptic.1} A cusp product-suspended operator
$A\in\Psi_{\cpsn{p}}^{m,m'}(Z;E,F)$ is said to be \emph{elliptic} if both
its symbol $\sigma(A)$ and its base family $\beta(A)$ are invertible; it is
\emph{fully elliptic} if its symbol $\sigma(A),$ its base family
$\beta(A)$ and its indicial family are all invertible.
\end{definition}


\providecommand{\bysame}{\leavevmode\hbox to3em{\hrulefill}\thinspace}
\providecommand{\MR}{\relax\ifhmode\unskip\space\fi MR }
\providecommand{\MRhref}[2]{%
  \href{http://www.ams.org/mathscinet-getitem?mr=#1}{#2}
}
\providecommand{\href}[2]{#2}


\end{document}